\documentclass[preprint, 11pt]{elsarticle}

\usepackage[T1]{fontenc}
\usepackage[utf8]{inputenc}
\usepackage[hmargin=0.9in]{geometry}
\usepackage[babel,german=quotes]{csquotes}
\usepackage{subfig}
\usepackage{graphicx}
\usepackage[colorlinks=true,citecolor=blue,urlcolor=blue]{hyperref}
\usepackage{caption}
\usepackage{amsmath}          % writing mathematical formulas
\usepackage{amsthm}           % writing mathematical theorems 
\usepackage{amssymb}          % writing mathematical symbols
\usepackage{bm, amsfonts}     % writing bold mathematical symbols
\usepackage{xcolor}
\usepackage{fixmath}
\usepackage{tikz}
\usepackage{bm}
\usepackage{makecell}

\newtheorem{thm}{Theorem}
\newdefinition{example}{Example}
\newproof{pf}{Proof}
\newcommand{\proofend}{\hfill $\Box$}

\newcommand {\dx} {\,{\rm d}{\mathbf x}}

\newcommand {\ds} {\,{\rm d}{\mathrm s}}

  \newcommand{\R}{\mathbb{R}}
  \newdefinition{rmk}{Remark}
  
  \newcommand{\pd}[2]{\frac{\partial #1}{\partial #2}}

  \newcommand{\td}[2]{\frac{\mathrm d #1}{\mathrm d #2}}

\newcommand{\beq}{\begin{equation}}
\newcommand{\eeq}{\end{equation}}

\makeatletter
\def\ps@pprintTitle{%
  \let\@oddhead\@empty
  \let\@evenhead\@empty
  \def\@oddfoot{
    \footnotesize\itshape
    \hfill\today
  }%
  \let\@evenfoot\@oddfoot}
\makeatother

\begin{document}

\begin{frontmatter} % statt \maketitle um email adressen anzuzeigen
  \title{Algebraic entropy fixes and convex limiting for continuous
 finite element discretizations of scalar hyperbolic conservation laws}

\author[TUDo]{Dmitri Kuzmin\corref{cor1}}
\ead{kuzmin@math.uni-dortmund.de}
\cortext[cor1]{Corresponding author}

\address[TUDo]{Institute of Applied Mathematics (LS III), TU Dortmund University\\ Vogelpothsweg 87,
  D-44227 Dortmund, Germany}

\author[KAUST]{Manuel Quezada de Luna}
\ead{manuel.quezada@kaust.edu.sa}

\address[KAUST]{King Abdullah University of Science and Technology (KAUST)\\ Thuwal 23955-6900, Saudi Arabia}

\journal{Computers and Mathematics with Applications}

\begin{abstract}
In this work, we modify a continuous Galerkin discretization of a scalar
hyperbolic conservation law using new algebraic correction
procedures. Discrete entropy conditions are used to determine the
minimal amount of entropy stabilization and constrain antidiffusive
corrections of a property-preserving low-order scheme. The addition
of a second-order entropy dissipative component to the antidiffusive
part of a nearly entropy conservative numerical flux is generally
insufficient to prevent violations of local bounds in shock regions.
Our monolithic convex limiting technique adjusts a given target flux
in a manner which guarantees preservation of invariant domains, validity
of local maximum principles, and entropy stability. The new methodology
combines the advantages of modern entropy stable / entropy conservative
schemes and their local extremum diminishing counterparts. The process of
algebraic flux correction is based on inequality constraints which
provably provide the desired properties. No free parameters are involved.
The proposed algebraic fixes are readily applicable to
unstructured meshes, finite element methods, general time
discretizations, and steady-state residuals. Numerical studies of
explicit entropy-constrained
schemes are performed for linear and nonlinear test problems.

\end{abstract}
\end{frontmatter}

\section{Introduction}

Entropy stability \cite{chen,ray,tadmor,tadmor2016} and
preservation of invariant domains  \cite{Guermond2018,Guermond2016,convex} play an important role in the design of numerical methods for nonlinear hyperbolic conservation laws. A failure to comply with these design criteria may result in nonphysical artefacts and/or convergence to wrong weak solutions. Modern high-resolution schemes are commonly equipped with flux or slope limiters which guarantee the validity of discrete maximum principles but may fail to satisfy entropy conditions. On the other hand, entropy stability of a high-order method does not guarantee the invariant domain preservation (IDP) property and numerical solutions may exhibit undershoots/overshoots.

Recent years have witnessed an increased interest of the finite element
community in analysis and design of algebraic flux correction (AFC) schemes
\cite{afc_analysis1,afc_analysis2,afc1,afc2}.
The AFC methodology modifies a standard Galerkin discretization by adding
artificial diffusion operators and limited antidiffusive fluxes. The convex
limiting techniques proposed in \cite{Guermond2018,convex,convex2} are applicable to nonlinear hyperbolic problems and lead to high-order IDP approximations.
However, additional inequality constraints must be taken into account to
ensure entropy stability. In the context of finite volume and discontinuous
Galerkin (DG) approximations, entropy stability is commonly achieved by
adding some entropy viscosity to an entropy conservative numerical flux.
For a comprehensive review of entropy stable schemes based on this
design philosophy, we refer to Tadmor \cite{tadmor,tadmor2016}. A
representation of continuous finite element approximations in terms
of numerical fluxes is also possible \cite{Selmin1993,Selmin1996}
but rather uncommon and requires the use of edge-based data structures
\cite{lohner_book}. Therefore, the use of formulations that add
diffusive fluxes to the residual of the Galerkin discretization is
preferred in the AFC literature \cite{afc1}.

As shown by Guermond et al. \cite{Guermond2014}, entropy stability is
an essential requirement for convergence of AFC schemes to correct weak
solutions of nonlinear hyperbolic problems.
Residual-based entropy viscosity \cite{Guermond2018,Guermond2014} was
found to be a good way to stabilize
flux-corrected continuous Galerkin (CG)
approximations \cite{Guermond2018,convex2}.
However, it involves a free parameter and does not guarantee entropy
stability. The entropy fixes proposed by
Abgrall et al. \cite{abgrall,ranocha} use Rusanov-type 
dissipation
terms to enforce a cell entropy inequality. In
contrast to finite volume and DG methods, construction of entropy
conservative CG schemes for which this inequality holds as equality is an open
problem. Hence, the minimal amount of entropy stabilization needs
to be determined without enforcing local entropy conservation.
The entropy stability
conditions that we use in the present paper are derived by
adapting Tadmor's \cite{tadmor} design criteria to
the CG setting. The key ingredients
of the proposed methodology are
\begin{itemize}
\item inequality constraints that guarantee entropy stability
  and preservation of invariant domains;
\item a general framework for designing algebraic flux correction
 schemes based on such constraints;
\item new parameter-free algorithms for construction 
  and limiting of entropy stabilization terms.
\end{itemize}
We begin with the CG space discretization of the initial value problem
in Section 2.
After introducing the new AFC tools and their theoretical foundations
in Sections 3-5, we summarize the proposed algorithm in Section 6,
perform numerical studies in Section 7, and draw conclusions in Section~8.

\section{Finite element discretization}

Let $u(\mathbf{x},t)$ be a scalar conserved
quantity depending on the space location $\mathbf{x}\in
\R^d,\ d\in\{1,2,3\}$ and time instant $t\ge 0$. Consider
an initial value problem of the form 
\begin{subequations}
\begin{align}
 \pd{u}{t}+\nabla\cdot\mathbf{f}(u)=0 &\qquad\mbox{in}\ \R^d\times\R_+,
\label{ibvp-pde}\\
 u(\cdot,0)=u_0 &\qquad\mbox{in}\ \R^d,\label{ibvp-ic}
\end{align}
\end{subequations}
where $\mathbf{f}=(\mathsf{f}_1,\ldots,\mathsf{f}_d)$ is a possibly
nonlinear flux function and $u_0:\R^d\to\mathcal G$ is an initial data
belonging to a convex set $\mathcal G$. The set $\mathcal G$ is
called an {\it invariant set} of problem  \eqref{ibvp-pde}--\eqref{ibvp-ic}
if the exact solution $u$ stays in $\mathcal G$ for all $t> 0$ \cite{Guermond2016}.
A convex function
$\eta:\mathcal G\to\R$ is called an {\it entropy} and $v=\eta'$ is called an
{\it entropy}
variable if there exists an {\it entropy flux} 
$\mathbf{q}:\mathcal G\to\mathbb{R}^d$ such that $v(u)\mathbf{f}'(u)=\mathbf{q}'(u)$.
 A weak solution
$u$ of \eqref{ibvp-pde} is called an {\it entropy solution}
if the entropy inequality
\beq
 \pd{\eta}{t}+\nabla\cdot\mathbf{q}(u)\le 0 \qquad\mbox{in}\ \R^d\times\R_+
\label{ent-ineq}
\eeq
holds for any {\it entropy pair} $(\eta,\mathbf{q})$. For any smooth weak
solution, the conservation law
\beq
 \pd{\eta}{t}+\nabla\cdot\mathbf{q}(u)=0 \qquad\mbox{in}\ \R^d\times\R_+
 \eeq
 can be derived from \eqref{ibvp-pde} using multiplication by the entropy
 variable $v$, the chain rule, and the definition of an entropy pair. Hence,
 entropy is conserved in smooth regions and dissipated at shocks.

Adopting the terminology of
Guermond et al. \cite{Guermond2018,Guermond2016}, we will call 
a numerical scheme {\it invariant
domain preserving} (IDP) if the solution of the (semi-)discrete problem
is guaranteed to stay in an invariant set $\mathcal G$. Additionally, a property-preserving
discretization of \eqref{ibvp-pde} should be {\it entropy stable}, i.e., it should satisfy
a discrete version of the entropy inequality \eqref{ent-ineq}. The lack of
entropy stability is a typical reason for convergence of numerical schemes
to nonphysical weak solutions.

Restricting the spatial domain to $\Omega\subset\R^d$ and imposing periodic
boundary conditions for simplicity,
we discretize \eqref{ibvp-pde} in space using a
conforming mesh $\mathcal T_h=\{K_1,\ldots,K_{E_h}\}$ of linear
($\mathbb{P}_1$) or multilinear ($\mathbb{Q}_1$)
finite elements. The globally continuous basis
functions $\varphi_1,\ldots,\varphi_{N_h}$
are associated with the vertices $\mathbf{x}_1,\ldots,\mathbf{x}_{N_h}$
of $\mathcal T_h$. Let $\mathcal E_i$
denote the set of (numbers of) elements containing the vertex
$\mathbf{x}_i$ and $\mathcal N^e$ is the set of (numbers of)
nodes belonging to $K^e$. The computational stencil of node $i$
is the integer set $\mathcal N_i=\bigcup_{e\in\mathcal E_i}\mathcal N^e$.
Substituting the
finite element approximations
\beq
u_h=\sum_{j=1}^{N_h}u_j\varphi_j,\qquad
\mathbf{f}_h=\sum_{j=1}^{N_h}\mathbf{f}_j\varphi_j\approx \mathbf{f}(u_h)
\eeq
into the weak form of \eqref{ibvp-pde} and using $\varphi_i,\ i\in
\{1,\ldots,N_h\}$ as a
test function, we obtain \cite{convex}
\beq\label{semi-high}
\sum_{e\in\mathcal E_i}\sum_{j\in\mathcal N^e}m_{ij}^e\td{u_j}{t}
=-
\sum_{e\in\mathcal E_i}\sum_{j\in\mathcal N^e}\mathbf{c}_{ij}^e\cdot\mathbf{f}_j
=-
\sum_{e\in\mathcal E_i}\sum_{j\in\mathcal N^e\backslash\{i\}}\mathbf{c}_{ij}^e\cdot
(\mathbf{f}_j-\mathbf{f}_i),
\eeq
\beq
m_{ij}^e=\int_{K^e}\varphi_i\varphi_j\dx,\qquad
\mathbf{c}_{ij}^e=\int_{K^e}\varphi_i\nabla\varphi_j\dx,\qquad
\sum_{j\in\mathcal N^e}\varphi_j(\mathbf{x})=1\ \forall
\mathbf{x}\in K^e.
\eeq
The choice of the time integration method should ensure at least conditional
$L^2$ stability of the fully discrete problem for linear flux functions of
the form $\mathbf{f}(u)=\mathbf{v}u,
\ \mathbf{v}\in\R^d$. The lack of nonlinear stability can be cured using
the algebraic flux correction tools that we present in the next sections.

\section{Property-preserving flux correction}\label{sec:prop_preserving}

To enforce entropy inequalities and local discrete maximum
principles, we
approximate \eqref{semi-high} by 
\beq\label{semi-afc}
\sum_{e\in\mathcal E_i}m_{i}^e\td{u_i}{t}=
\sum_{e\in\mathcal E_i}\sum_{j\in\mathcal N^e\backslash\{i\}}
   [g_{ij}^e-
      \mathbf{c}_{ij}^e\cdot(\mathbf{f}_j-\mathbf{f}_i)],
\eeq
where
\beq
m_{i}^e=\sum_{j\in\mathcal N^e}m_{ij}^e=\int_{K^e}\varphi_i\dx
\eeq
are the diagonal entries of the lumped element mass matrix
and $g_{ij}^e$ are numerical fluxes such that
\beq
g_{ji}^e=-g_{ij}^e\qquad \forall i\in\mathcal N^e,\ j\in\mathcal N^e\backslash\{i\}.
\eeq
The standard continuous Galerkin scheme \eqref{semi-high} can be
written in the form \eqref{semi-afc} using the fluxes
\beq\label{gal-flux}
g_{ij}^{e,\rm CG}=m_{ij}^e(\dot u_i-\dot u_j).
\eeq
The nodal time derivatives $\dot u_i=\td{u_i}{t}$ are defined by
\eqref{semi-high}. To avoid inversion of the consistent mass 
matrix, an approximate solution of this linear system for
$\dot u$ can be obtained efficiently using a few Richardson's
iterations preconditioned by the lumped mass matrix \cite{donea,quartapelle}. This
approach to calculating $\dot u$ corresponds to an approximation by a
truncated Neumann series \cite{Guermond2014,quartapelle}.

The purpose of algebraic flux correction (AFC) is to replace $g_{ij}^{e,\rm CG}$
with a flux that contains enough numerical dissipation to ensure
preservation of invariant domains, validity of local discrete maximum
principles and/or entropy stability. On the other hand,
the levels of numerical diffusion should be kept small enough to achieve
optimal convergence behavior for problems with smooth exact solutions.
Similarly to PDE-constrained optimization problems, AFC schemes are
designed to adjust the {\it control variables} $g_{ij}^{e}$ in a
way which guarantees the validity of certain constraints for the
{\it state variables} $u_i$ while staying as close as possible
to a given {\it target}. Numerical solution of global
constrained optimization problems is feasible \cite{bochev}
but costly. Therefore, we will design $g_{ij}^{e}$ using sufficient
conditions (box constraints) to derive simple closed-form
approximations which provide the desired properties.
\medskip

Let $(\eta,\mathbf{q})$ be an entropy pair and $v=\eta'(u)$
the corresponding entropy variable. Define
\beq
\boldsymbol{\psi}(u)=v(u)\mathbf{f}(u)-\mathbf{q}(u).
\eeq
A sufficient condition for entropy stability of the
semi-discrete problem \eqref{semi-afc} is given by (cf. \cite{chen,pazner})
\beq\label{condES}
\frac{v_i-v_j}{2}[g_{ij}^{e}-\mathbf{c}_{ij}^e\cdot(\mathbf{f}_j+\mathbf{f}_i)]
\le \mathbf{c}_{ij}^e\cdot[\boldsymbol{\psi}(u_j)
  -\boldsymbol{\psi}(u_i)].
\eeq
In the following Theorem, we show that \eqref{condES} implies
the validity of a semi-discrete entropy inequality.

\begin{thm}[Entropy stability of AFC schemes
    \cite{Guermond2018,Guermond2016}]\label{thm1}
If condition \eqref{condES} holds for each flux $g_{ij}^e$, then the solution 
of the semi-discrete problem \eqref{semi-afc} satisfies the discrete entropy
inequality
\beq
\sum_{e\in\mathcal E_i}m_{i}^e\td{\eta(u_i)}{t}\le
\sum_{e\in\mathcal E_i}\sum_{j\in\mathcal N^e\backslash\{i\}}
   [G_{ij}^e-
      \mathbf{c}_{ij}^e\cdot(\mathbf{q}_j-\mathbf{q}_i)],
\eeq
where 
\beq
G_{ij}^e=\frac{v_i+v_j}{2}\,g_{ij}^{e}-
  \frac{v_i-v_j}{2}\,\mathbf{c}_{ij}^e\cdot (\mathbf{f}_j-\mathbf{f}_i).
\eeq
\end{thm} 

\begin{pf}
  We have  $\sum_{e\in\mathcal E_i}m_{i}^e\td{\eta(u_i)}{t}
  =\sum_{e\in\mathcal E_i}m_{i}^ev_i\td{u_i}{t}
  =\sum_{e\in\mathcal E_i}v_i\sum_{j\in\mathcal N^e\backslash\{i\}}[g_{ij}^{e}
    -\mathbf{c}_{ij}^e\cdot(\mathbf{f}_j-\mathbf{f}_i)]$
  by the chain rule and the definition of the entropy variable
  $v_i=\eta'(u_i)$. Using
  the zero sum property
   \mbox{$\sum_{j\in\mathcal N^e}\mathbf{c}_{ij}^e=\mathbf{0}$}
   of the discrete gradient operator, and the stability
   condition  \eqref{condES},  we find that 
\begin{align*}
&  v_i\sum_{j\in\mathcal N^e\backslash\{i\}}[g_{ij}^{e}
    -\mathbf{c}_{ij}^e\cdot(\mathbf{f}_j-\mathbf{f}_i)]
  = v_i\sum_{j\in\mathcal N^e\backslash\{i\}}[g_{ij}^{e}-\mathbf{c}_{ij}^e\cdot(\mathbf{f}_j+\mathbf{f}_i)]-2v_i\mathbf{c}_{ii}^e\cdot\mathbf{f}_i\\
  &= \sum_{j\in\mathcal N^e\backslash\{i\}}\left(
\frac{v_i+v_j}{2}
     [g_{ij}^{e}-\mathbf{c}_{ij}^e\cdot(\mathbf{f}_j+\mathbf{f}_i)]+\frac{v_i-v_j}{2}
     [g_{ij}^{e}-\mathbf{c}_{ij}^e\cdot(\mathbf{f}_j+\mathbf{f}_i)]\right)
  -2v_i\mathbf{c}_{ii}^e\cdot\mathbf{f}_i\\
  &\le \sum_{j\in\mathcal N^e\backslash\{i\}}\left(\frac{v_i+v_j}{2}[g_{ij}^{e}
    -\mathbf{c}_{ij}^e\cdot(\mathbf{f}_j+\mathbf{f}_i)]
+\mathbf{c}_{ij}^e\cdot
[\boldsymbol{\psi}(u_j)-\boldsymbol{\psi}(u_i)]
\right)-2v_i\mathbf{c}_{ii}^e\cdot\mathbf{f}_i\\
&=\sum_{j\in\mathcal N^e\backslash\{i\}}\left(\frac{v_i+v_j}{2}[g_{ij}^{e}
  -\mathbf{c}_{ij}^e\cdot(\mathbf{f}_j+\mathbf{f}_i)]
+\mathbf{c}_{ij}^e\cdot
[\boldsymbol{\psi}(u_j)+\boldsymbol{\psi}(u_i)]\right) -2\mathbf{c}_{ii}^e\cdot[
v_i\mathbf{f}_i-\boldsymbol{\psi}(u_i)]\\
(*)\quad &=\sum_{j\in\mathcal N^e\backslash\{i\}}\left(
\frac{v_i+v_j}{2}g_{ij}^{e}
- \mathbf{c}_{ij}^e\cdot\left[
  \frac{v_i-v_j}{2}(\mathbf{f}_j-\mathbf{f}_i)+
 \mathbf{q}_j+\mathbf{q}_i\right]\right)
-2\mathbf{c}_{ii}^e\cdot\mathbf{q}_i\\
&=\sum_{j\in\mathcal N^e\backslash\{i\}}\left(
\frac{v_i+v_j}{2}g_{ij}^{e}- \mathbf{c}_{ij}^e\cdot\left[
  \frac{v_i-v_j}{2}(\mathbf{f}_j-\mathbf{f}_i)+
  \mathbf{q}_j-\mathbf{q}_i\right]\right)=
\sum_{j\in\mathcal N^e\backslash\{i\}}
   [G_{ij}^e-
      \mathbf{c}_{ij}^e\cdot(\mathbf{q}_j-\mathbf{q}_i)].
\end{align*}
Summing over $e\in\mathcal E_i$, we conclude that the assertion of
the Theorem is true.
    \qquad\proofend
\end{pf}

By definition of $\mathbf{c}_{ij}^e$, we have
$\sum_{e\in\mathcal E_i}\mathbf{c}_{ij}^e=-\sum_{e\in\mathcal E_i}\mathbf{c}_{ji}^e$ if $i$
or $j$ is an interior node. Under the assumption of periodic boundary conditions,
this property holds for all nodes. In particular, we have $\sum_{e\in\mathcal E_i}\mathbf{c}_{ii}^e=\mathbf{0}$. Using the identity marked by (*) in the proof of Theorem \ref{thm2}, we obtain the estimate
\begin{align}\label{detadt}
\sum_{i=1}^{N_h}\sum_{e\in\mathcal E_i}
 m_i^e\td{\eta_i}{t}&\le\sum_{i=1}^{N_h}\sum_{j=1\atop j>i}^{N_h}
   \frac{v_i+v_j}{2}\sum_{e\in\mathcal E_i}\underbrace{(g_{ij}^{e}+g_{ji}^{e})}_{=0}
   \\
&-\sum_{i=1}^{N_h}
\sum_{j=1\atop j>i}^{N_h}
\left[\frac{v_i-v_j}{2}(\mathbf{f}_j-\mathbf{f}_i)+
  \mathbf{q}_j+\mathbf{q}_i\right]\cdot
\underbrace{\sum_{e\in\mathcal E_i}(\mathbf{c}_{ij}^e+\mathbf{c}_{ji}^e)}_{=0}
-2\mathbf{q}_i\cdot\underbrace{\sum_{e\in\mathcal E_i}\mathbf{c}_{ii}^e}_{=\mathbf{0}}
=0\nonumber
\end{align}
in accordance with the fact that $\td{}{t}\int_{\Omega}\eta(u)\dx
\le 0$ for the entropy solution of an initial boundary-value problem
with $\int_{\partial\Omega}\mathbf{q}(u)\cdot\mathbf{n}\ds=0$, where $\mathbf{n}$
denotes the unit outward normal. Note that the validity of estimate
\eqref{detadt}
for the square entropy $\eta=\frac{u^2}{2}$ implies $L^2$ stability.

Suppose that the exact entropy solution $u$ belongs to a convex invariant
set $\mathcal
G=[u^{\min},u^{\max}]$. Then
a semi-discrete scheme of the form \eqref{semi-afc} is
invariant domain preserving (IDP) if it satisfies
\beq
u^{\min}\le u_i^{\min}(t)\le u_i(t)\le u_i^{\max}(t)\le u^{\max}
\qquad \forall t\ge 0.
\eeq
A fully discrete scheme possesses the IDP property if similar inequality
constraints hold at each discrete time level $t^{n}=n\Delta t,\ n\in\mathbb{N}$ or
stage of a strong stability preserving (SSP) Runge-Kutta method \cite{ssprev}.
The following Theorem provides a sufficient condition for
the design of IDP approximations.

\begin{thm}[Guermond-Popov IDP criterion
    \cite{Guermond2018,Guermond2016}]\label{thm2}
  Consider a semi-discrete scheme of the form
\beq\label{bar-idp}
\sum_{e\in\mathcal E_i}m_{i}^e\td{u_i}{t}=
\sum_{e\in\mathcal E_i}\sum_{j\in\mathcal N^e\backslash\{i\}}
  2d_{ij}^e(\bar u_{ij}^e-u_i),\qquad i\in\{1,\ldots,N_h\},
  \eeq
  where $m_i^e>0$ and
  $d_{ij}^e>0$ for all $j\in\mathcal N^e\backslash\{i\}$. Let $\mathcal G$
  be a convex set. Assume that 
  \beq 
  u_i\in\mathcal G,\qquad \bar u_{ij}^e\in\mathcal G\quad
  \forall j\in\mathcal N^e\backslash\{i\}.
  \eeq
If the time step $\Delta t$ satisfies
   \beq\label{bar-cfl}
\Delta t\sum_{e\in\mathcal E_i}
\sum_{j\in\mathcal N^e\backslash\{i\}}2d_{ij}^e\le m_i=
\sum_{e\in\mathcal E_i}m_i^e,
  \eeq
 then an explicit SSP Runge-Kutta time discretization
  of \eqref{bar-idp} is IDP w.r.t. $\mathcal G$. 
\end{thm}

\begin{pf}
  Each stage of an explicit SSP-RK method is a forward Euler update
  of the form 
 \begin{align*}
  \bar u_i&=u_i+\frac{\Delta t}{m_i}
\sum_{e\in\mathcal E_i}\sum_{j\in\mathcal N^e\backslash\{i\}}
2d_{ij}^e(\bar u_{ij}^e-u_i)\\
& =\left(1-\frac{\Delta t}{m_i}\sum_{e\in\mathcal E_i}
\sum_{j\in\mathcal N^e\backslash\{i\}}2d_{ij}^e\right)u_i
+\frac{\Delta t}{m_i}\sum_{e\in\mathcal E_i}\sum_{j\in\mathcal N^e\backslash\{i\}}
2d_{ij}^e\bar u_{ij}^e.
 \end{align*}
 Under the time step restriction \eqref{bar-cfl}, this 
 representation of
 the fully discrete scheme 
 implies that $\bar u_i$ is a convex combination of $u_i\in\mathcal G$ and
$\bar u_{ij}^e\in\mathcal G$. Since $\mathcal G$ is convex, the result
 $\bar u_i$ stays in  $\mathcal G$ \cite{Guermond2016}.

    \qquad\proofend
\end{pf}

\begin{rmk} After the global residual  assembly, the semi-discrete AFC scheme
 \eqref{semi-afc} becomes
  \beq
m_{i}\td{u_i}{t}=\sum_{j\in\mathcal N_i}
   [g_{ij}-
      \mathbf{c}_{ij}\cdot(\mathbf{f}_j-\mathbf{f}_i)],
\eeq
where
\beq
m_{i}=\sum_{j\in\mathcal N_i}m_{ij},\qquad
m_{ij}=\sum_{e\in\mathcal E_i\cap\mathcal E_j}m_{ij}^e,\qquad
\mathbf{c}_{ij}=\sum_{e\in\mathcal E_i\cap\mathcal E_j}\mathbf{c}_{ij}^e.
\eeq
In this paper, we assemble $g_{ij}=\sum_{e\in\mathcal E_i}g_{ij}^e$
from element contributions $g_{ij}^e$. However, the corrected flux can
also be determined directly.
In AFC schemes of this kind, inequality constraints for $g_{ij}$
are formulated using the coefficients of global matrices
(cf. \cite{Guermond2018,afc1,convex,fctools}). All
algorithms to be presented below
can be easily converted to the post-assembly format by
dropping the superscript~$e$.
\end{rmk}

\section{Entropy stable AFC schemes}\label{sec:ent_stable_afc}

Let us begin with the derivation of semi-discrete AFC schemes
satisfying the entropy stability condition \eqref{condES}.
A low-order approximation of local Lax-Friedrichs (LLF)
type is defined by 
\beq\label{glow}
g_{ij}^{e,\rm LLF}=d_{ij}^{e,\max}(u_j-u_i),
\eeq
where
\beq\label{dijmax}
d_{ij}^{e,\max}=\begin{cases}
\max\{|\mathbf{c}_{ij}^e|,
|\mathbf{c}_{ji}^e|\}\max\{\lambda_{ij}^{\max},\lambda_{ji}^{\max}\} & \mbox{if}\
i\in\mathcal N^e,\ j\in \mathcal N^e\backslash\{i\},\\
-\sum_{k\in \mathcal N^e\backslash\{i\}}d_{ik}^e & \mbox{if}\ j=i\in\mathcal N^e,\\
0 & \mbox{otherwise}
\end{cases}
\eeq
are artificial diffusion coefficients proportional to the
maximum wave speed \cite{Guermond2016,convex,convex2}
\beq
\lambda_{ij}^{\max}=
\max_{\omega\in[0,1]}|\mathbf{n}_{ij}^e\cdot
\mathbf{f}'(\omega u_i+(1-\omega) u_j)|,\qquad \mathbf{n}_{ij}^e=
\frac{\mathbf{c}_{ij}^e}{|\mathbf{c}_{ij}^e|}.
\eeq
The LLF flux defined by \eqref{glow} satisfies \eqref{condES},
as shown by Chen and Shu \cite{chen} in the context of
entropy stable DG methods. Guermond and Popov \cite{Guermond2016}
proved that the fully discrete SSP-RK version of the LLF-AFC scheme
is IDP and satisfies a discrete entropy inequality for any
entropy pair $(\eta,\mathbf{q})$. However, the accuracy of
the LLF approximation is first-order at best and significant
amounts of numerical diffusion can be removed without losing
the entropy stability property.

To derive a flux control $g_{ij}^e$ corresponding to a second-order entropy
stable approximation, we adopt Tadmor's \cite{tadmor87,tadmor} design
philosophy which is based on comparison with entropy conservative schemes.
Suppose that condition \eqref{condES} holds as equality for some $g_{ij}^{e,\rm EC}$.
Then it holds as inequality for
\beq
g_{ij}^{e,\rm ES}=g_{ij}^{e,\rm EC}+\nu_{ij}^e(v_j-v_i),
\eeq
where $\nu_{ij}^e\ge 0$ is an entropy viscosity coefficient. This simple comparison
principle provides a powerful tool for the design of entropy stable finite volume
\cite{fjordholm,ray,tadmor2016} and DG \cite{chen,pazner} methods. For our AFC
scheme \eqref{semi-afc} to be entropy conservative, the fluxes
 $g_{ij}^{e,\rm EC}=-g_{ji}^{e,\rm EC}$ would need to satisfy 
\begin{align*}
  \frac{v_i-v_j}{2}[g_{ij}^{e,\rm EC}-\mathbf{c}_{ij}^e\cdot(\mathbf{f}_j+\mathbf{f}_i)]
&=\mathbf{c}_{ij}^e\cdot[\boldsymbol{\psi}(u_j)
  -\boldsymbol{\psi}(u_i)],\\
\frac{v_i-v_j}{2}[g_{ij}^{e,\rm EC}+\mathbf{c}_{ji}^e\cdot(\mathbf{f}_j+\mathbf{f}_i)]
&=\mathbf{c}_{ji}^e\cdot[\boldsymbol{\psi}(u_i)
  -\boldsymbol{\psi}(u_j)].
\end{align*}
Since this system is overdetermined in the case $\mathbf{c}_{ij}^e\ne -\mathbf{c}_{ji}^e$,
we perform element-level flux correction using generalized entropy-stable target
fluxes of the form
\beq\label{gstab}
g_{ij}^{e,\rm ES}=d_{ij}^{e,\min}(u_j-u_i)+\nu_{ij}^e(v_j-v_i),
\eeq
where $d_{ij}^{e,\min}\in [0,d_{ij}^{e,\max}]$ is the minimal nonnegative diffusion coefficient satisfying
the symmetry condition $d_{ij}^{e,\min}=d_{ji}^{e,\min}$ and 
condition \eqref{condES} for both nodes. The value of $d_{ij}^{e,\min}$ is given by
\beq\label{dijdef}
d_{ij}^{e,\min}=\begin{cases}
\frac{\min\{Q_{ij}^e,0,Q_{ji}^e\}}{(v_i-v_j)(u_j-u_i)} &\mbox{if}\ u_i\ne u_j,\\
0 & \mbox{if}\ u_i=u_j,
\end{cases}
\eeq
where
\beq\label{qijdef}
Q_{ij}^e=2\mathbf{c}_{ij}^e\cdot \left[\boldsymbol{\psi}(u_j)-\boldsymbol{\psi}(u_i)
+(v_i-v_j)\frac{\mathbf{f}_j+\mathbf{f}_i}{2}\right].
\eeq
By the mean value theorem, we have $v_i-v_j=\eta'(u_i)-\eta'(u_j)
=\eta''(\xi)(u_i-u_j)$  for some $\xi\in\R$. It follows that
$(v_i-v_j)(u_j-u_i)\le 0$ and, therefore, $d_{ij}^{e,\min}\ge 0$
for any convex entropy $\eta$.

\begin{rmk}
  In the absence of rounding errors, we have $d_{ij}^{e,\min}\le d_{ij}^{e,\max}$ by
  definition. In practice, division by a small number $(v_i-v_j)(u_j-u_i)$ may
  produce $d_{ij}^{e,\min}>d_{ij}^{e,\max}$ in regions where the numerical solution
  is almost constant. To avoid this, it is worthwhile to use $d_{ij}^{e,\max}$ as
  upper bound for $d_{ij}^{e,\min}$ in practical implementations. We also
  remark that the direct calculation of the diffusive flux $d_{ij}^{e,\min}(u_j-u_i)
  =\frac{\min\{Q_{ij}^e,0,Q_{ji}^e\}}{v_i-v_j}$ is less sensitive
  to rounding errors in the limit $|u_i-u_j|\to 0$.
  \end{rmk}

The AFC scheme corresponding to \eqref{gstab} with $\nu_{ij}^e=0$ is {\it barely entropy stable}.
Building on Tadmor's \cite{tadmor} ideas, we define the additional flux $\nu_{ij}^e(v_j-v_i)$
using the entropy viscosity coefficient
\beq\label{nuijdef}
\nu_{ij}^e=\max\left\{\frac{\mathbf{c}_{ij}^e\cdot
  \left[\mathbf{f}_j+\mathbf{f}_i
    -2\mathbf{f}\left(\frac{u_j+u_i}{2}\right)\right]}{v_j-v_i},0,
\frac{\mathbf{c}_{ji}^e\cdot\left[\mathbf{f}_j+\mathbf{f}_i
    -2\mathbf{f}\left(\frac{u_j+u_i}{2}\right)\right]}{v_i-v_j}
  \right\}
\eeq
which vanishes for linear flux functions $\mathbf{f}(u)$ and preserves
second-order accuracy for nonlinear ones.

\begin{rmk}
  For reasons explained in Remark 2, we recommend direct
calculation of the flux
  \beq
\nu_{ij}^e(v_j-v_i)=S_{ij}\max\left\{S_{ij}\mathbf{c}_{ij}^e\cdot
2\Delta \mathbf{f}_{ij},0,-S_{ij}
\mathbf{c}_{ji}^e\cdot 2\Delta \mathbf{f}_{ij} \right\},
  \eeq
  where $S_{ij}$ is the sign of $v_j-v_i$ and $\Delta \mathbf{f}_{ij}
  =\frac12(\mathbf{f}_j+\mathbf{f}_i)
    -\mathbf{f}\left(\frac{u_j+u_i}{2}\right)$ is the flux difference. 
  \end{rmk}

The use of \eqref{gstab} with $d_{ij}^{e,\min}$ defined by \eqref{dijdef} and $\nu_{ij}^e$
defined by \eqref{nuijdef} yields an entropy stable approximation which exhibits 
the desired convergence behavior for smooth data but may produce undershoots and/or
overshoots in the neighborhood of shocks. To enforce the IDP property and preservation
of local bounds, we use the monolithic
convex limiting techniques presented in the next section.

\section{Bound-preserving AFC schemes}

The highly dissipative LLF flux \eqref{glow} satisfies not only the entropy
condition \eqref{condES} but also the assumptions of Theorem \ref{thm2}. The corresponding
IDP bar states are given by
\beq\label{barstate}
\bar u_{ij}^e
=\frac{u_j+u_i}{2}-\frac{\mathbf{c}_{ij}^e
        \cdot(\mathbf{f}_j-\mathbf{f}_i)}{2d_{ij}^{e,\max}},
\eeq
where $d_{ij}^{e,\max}$ is defined by \eqref{dijmax}. Using the mean value theorem, one can show that \cite{convex}
\beq\label{bar_bounds}
\min\{u_i,u_j\}\le \bar u_{ij}^e
\le \max\{u_i,u_j\}.
\eeq
Hence, the algebraic LLF scheme is IDP by Theorem 2. To limit the raw antidiffusive
part
\beq
f_{ij}^e=g_{ij}^e-g_{ij}^{e,\rm LLF}
\eeq
of a given flux control $g_{ij}^e$ in a manner which preserves the IDP property, we
define
\beq\label{gidp}
g_{ij}^{e,*}=g_{ij}^{e,\rm LLF}+f_{ij}^{e,*}=d_{ij}^{e,\max}(u_j-u_i)+f_{ij}^{e,*}
\eeq
using the inequality-constrained antidiffusive flux \cite{convex,convex2}
\beq\label{fij_lim}
f_{ij}^{e,*}=\begin{cases}
  \min\,\left\{f_{ij}^e,2d_{ij}^e\min\,\{u_i^{\max}-\bar u_{ij}^e,
    \bar u_{ji}^e-u_j^{\min}\}\right\} & \mbox{if}\  f_{ij}^e>0,\\[0.25cm]
    \max\left\{f_{ij}^e,2d_{ij}^e\max\{u_i^{\min}-\bar u_{ij}^e,
    \bar u_{ji}^e-u_j^{\max}\}\right\} & \mbox{otherwise}.
\end{cases}
\eeq
This monolithic convex limiting (MCL) strategy was proposed in \cite{convex}. It guarantees that
\beq\label{afc-bounds}
\min_{j\in \mathcal N_i}u_j=:u_i^{\min}
\le \bar u_{ij}^{e,*}=\bar u_{ij}^e+\frac{f_{ij}^*}{2d_{ij}^{e,\max}}
\le u_i^{\max}:=\max_{j\in \mathcal N_i}u_j.
\eeq
The IDP property of the flux-corrected scheme can be shown using Theorem \ref{thm2},
see \cite{convex} for details.

\begin{rmk}
  A linearity-preserving version of the bounds $u_i^{\min}$
  and $u_i^{\max}$  can be constructed as proposed
  in Section 6.1 of \cite{convex}. The use of limiters that guarantee
  linearity preservation (i.e., produce  $f_{ij}^{e,*}=f_{ij}^e$
  for locally linear functions $u_h$)
  is essential for achieving optimal convergence to smooth
  solutions \cite{afc_lp}.
\end{rmk}

Formula \eqref{fij_lim} will leave the raw antidiffusive
flux $f_{ij}^e$ unchanged if it does not violate the AFC
inequality constraints \eqref{afc-bounds}. Hence, the quality
of flux-corrected solutions depends on the properties of the 
(stabilized) high-order method defined by \eqref{semi-afc}
with $g_{ij}^e=g_{ij}^{e,\rm LLF}+f_{ij}^e$. The target flux 
\beq\label{target-lumped}
f_{ij}^e=(d_{ij}^{e,\min}
-d_{ij}^{e,\max})(u_j-u_i)
+\nu_{ij}^e(v_j-v_i)
\eeq
corresponds to an entropy stable lumped-mass approximation. The addition
of \eqref{target-lumped} to $g_{ij}^{e,\rm LLF}$ replaces it with
$g_{ij}^{e,\rm ES}$. If limiting is performed using  \eqref{fij_lim},
the addition of $f_{ij}^{e,*}$ replaces $g_{ij}^{e,\rm LLF}$ with 
\beq
g_{ij}^{e,*}=(1-\alpha_{ij}^e)d_{ij}^{e,\max}(u_j-u_i)+\alpha_{ij}g_{ij}^{e,\rm ES},
\qquad \alpha_{ij}^e=\begin{cases}
\frac{f_{ij}^{e,*}}{f_{ij}^e} &\mbox{if}\ f_{ij}^e\ne 0,\\
0 & \mbox{otherwise}.
\end{cases}
\eeq
Recall that
$(v_i-v_j)(u_j-u_i)\le 0$ for any convex entropy $\eta$
by the mean value theorem and definition of the entropy
variable $v=\eta'(u)$. Since 
$\alpha_{ij}^e\in[0,1]$ and $d_{ij}^{e,\max}\ge d_{ij}^{e,\min}$, we have
\begin{align*}
(v_i-v_j)g_{ij}^{e,*}&=(v_i-v_j)[
(1-\alpha_{ij}^e)d_{ij}^{e,\max}(u_j-u_i)+
    \alpha_{ij}^ed_{ij}^{e,\min}(u_j-u_i)+ \alpha_{ij}^e\nu_{ij}^e(v_j-v_i)]\\
  & \le (v_i-v_j)d_{ij}^{e,\min}(u_j-u_i).
  \end{align*}
Thus the replacement of  $g_{ij}^{e,\rm LLF}$ by $g_{ij}^{e,*}$ produces an
entropy stable and bound-preserving approximation.

In the consistent-mass version of our AFC scheme, the target flux
to be used in \eqref{fij_lim} is given by
\beq\label{target}
f_{ij}^e=m_{ij}^e(\dot u_i-\dot u_j)+(d_{ij}^{e,\min}-d_{ij}^{e,\max})(u_j-u_i)+\nu_{ij}^e(v_j-v_i)
\eeq
and represents the antidiffusive part of
$g_{ij}^e=g_{ij}^{e,\rm CG}+g_{ij}^{e,\rm ES}$. The nodal
time derivatives $\dot u$ can be defined using \eqref{semi-high} or \eqref{semi-afc}.
In the numerical examples of Section \ref{sec:num}, we use the LLF approximation
\beq\label{udotdef}
\dot u_i=\frac{1}{m_i}
\sum_{e\in\mathcal E_i}\sum_{j\in\mathcal N^e\backslash\{i\}}
   [d_{ij}^{e,\max}(u_j-u_i)-\mathbf{c}_{ij}^e\cdot(\mathbf{f}_j-\mathbf{f}_i)].
\eeq
In addition to being rather inexpensive, it has the positive effect of introducing high-order background stabilization \cite{convex}
even for linear advection problems, for which
$\nu_{ij}^e$ defined by \eqref{nuijdef} vanishes.

The inclusion of $m_{ij}^e(\dot u_i-\dot u_j)$ may
require additional limiting of $f_{ij}^{e,*}$ to ensure that the final flux
\beq\label{gfin}
g_{ij}^{e,**}=g_{ij}^{e,\rm LLF}+f_{ij}^{e,**}=d_{ij}^{e,\max}(u_j-u_i)+f_{ij}^{e,**}
\eeq
of our property-preserving AFC scheme \eqref{semi-afc} will satisfy the entropy stability condition
\beq\label{condES2}
\frac{v_i-v_j}{2}[g_{ij}^{e,**}-\mathbf{c}_{ij}^e\cdot(\mathbf{f}_j+\mathbf{f}_i)]
\le \mathbf{c}_{ij}^e\cdot[\boldsymbol{\psi}(u_j)
  -\boldsymbol{\psi}(u_i)].
\eeq
Substituting \eqref{gfin} into \eqref{condES2}, we obtain a limiting
criterion for the algebraic entropy fix
\beq\label{fij_lim_fix}
f_{ij}^{e,**}
=\begin{cases}
\frac{\min\{Q_{ij}^{e,*},(v_i-v_j)f_{ij}^{e,*},Q_{ji}^{e,*}\}}{v_i-v_j}
&\mbox{if}\ (v_i-v_j)f_{ij}^{e,*}>0,\\
 f_{ij}^* & \mbox{otherwise},
\end{cases}
\eeq
where
\beq
Q_{ij}^{e,*}=Q_{ij}^{e}-(v_j-v_i)d_{ij}^{e,\max}(u_i-u_j)
\eeq
are upper bounds for entropy-producing  fluxes. The value of $Q_{ij}^{e}$
is given by \eqref{qijdef}. Note that
$$Q_{ij}^{e,*}
\ge (v_i-v_j)(d_{ij}^{e,\max}-d_{ij}^{e,\min})(u_i-u_j)=\eta''(\xi)(d_{ij}^{e,\max}-d_{ij}^{e,\min})
(u_i-u_j)^2$$ for some $\xi\in\R$. Since $d_{ij}^{e,\max}\ge d_{ij}^{e,\min}$, the
bounds $Q_{ij}^{e,*}$ are nonnegative for any convex entropy $\eta$.

\begin{rmk}
  The replacement of $f_{ij}^{e,*}$ by the entropy-corrected flux
  $f_{ij}^{e,**}$ is equivalent to multiplication by a
  correction factor $\alpha_{ij}^e\in[0,1]$. Hence, it does not affect the IDP property
  of our AFC scheme.
\end{rmk}

\section{Summary of the algorithm}
\label{sec:summ}

Let us now summarize the algorithmic steps to be performed and the properties of the AFC scheme that ensures entropy stability and preservation of local bounds. At the semi-discrete level, the flux-corrected CG approximation is defined by the nonlinear system of ordinary differential equations
\beq\label{semi-afc-fix}
m_i\td{u_i}{t}=
\sum_{e\in\mathcal E_i}\sum_{j\in\mathcal N^e\backslash\{i\}}
   [d_{ij}^{e,\max}(u_j-u_i)+f_{ij}^{e,**}-
     \mathbf{c}_{ij}^e\cdot(\mathbf{f}_j-\mathbf{f}_i)],\qquad
   i=1,\ldots,N_h,
   \eeq
   where $d_{ij}^{e,\max}$ is the maximal speed diffusion coefficient defined by \eqref{dijmax}. If all corrections are included, the computation of the limited antidiffusive flux $f_{ij}^{e,**}$ involves the following steps:
   \begin{enumerate}
   \item Calculate the minimal diffusion coefficient $d_{ij}^{e,\min}$ using \eqref{dijdef}.
   \item Calculate the entropy viscosity coefficient $\nu_{ij}^{e}$ using \eqref{nuijdef}.
   \item Calculate the approximate time derivatives $\dot u_i$ using \eqref{udotdef}.
   \item Calculate the raw antidiffusive fluxes $f_{ij}^e$ using \eqref{target}.  
   \item Calculate the local bounds $u_i^{\max}$ and $u_i^{\min}$ using
     \eqref{afc-bounds}.
   \item Calculate the bound-preserving fluxes $f_{ij}^{e,*}$ using \eqref{fij_lim}.
   \item Calculate the entropy-corrected fluxes $f_{ij}^{e,**}$ using \eqref{fij_lim_fix}.  
   \end{enumerate}
   The following implication of Theorem \ref{thm2} provides a sufficient condition for an explicit time discretization of the nonlinear semi-discrete AFC problem \eqref{semi-afc-fix} to be locally bound preserving.

   \begin{thm}[IDP property of the flux-corrected CG scheme]\label{thm3}
 An explicit SSP Runge-Kutta time discretization of system \eqref{semi-afc-fix} 
 satisfies the local maximum principle
 \beq
u_i^{\min}\le \bar u_i\le u_i^{\max}
 \eeq
 under the time step restriction
 \beq
\Delta t\sum_{e\in\mathcal E_i}
\sum_{j\in\mathcal N^e\backslash\{i\}}2d_{ij}^{e,\max}\le m_i.
  \eeq
\end{thm}

\begin{pf}
 To apply Theorem \ref{thm2},
  we notice that each SSP Runge-Kutta stage can be written as (cf. \cite{convex})
  \beq
  \bar u_i=u_i+\frac{\Delta t}{m_i}\sum_{e\in\mathcal E_i}
  \sum_{j\in\mathcal N^e\backslash\{i\}}2d_{ij}^{e,\max}(\bar
  u_{ij}^{e,**}-u_i),
  \eeq
  where
  \beq
  u_i^{\min}\le  \bar u_{ij}^{e,**}=\bar u_{ij}^e+\frac{f_{ij}^{e,**}}{2d_{ij}^{e,\max}}
  \le  u_i^{\max}
  \eeq
  by virtue of \eqref{bar_bounds}, \eqref{fij_lim}, and \eqref{fij_lim_fix}. The
  desired result follows by the convexity argument.
    \qquad\proofend
\end{pf}

\begin{rmk}
  The applicability of the presented AFC tools is not restricted to explicit
  SSP Runge-Kutta time discretizations. 
Nonlinear discrete problems associated with implicit time discretizations
and the steady state limit of \eqref{semi-afc-fix} 
can be analyzed as in \cite{afc_analysis1,CL-diss}, see also the Appendix of \cite{convex}.
  \end{rmk}

\section{Numerical examples}\label{sec:num}
In this section, we perform numerical experiments for linear and nonlinear scalar problems.
The purpose of this numerical study is to demonstrate that the proposed methodology provides optimal accuracy for linear ($\mathbb{P}_1$) and multilinear ($\mathbb{Q}_1$) finite element approximations to smooth solutions, and that it behaves as expected on structured as well as unstructured meshes. In the description of the numerical results, we use the abbreviation LO-ES-IDP for the low-order LLF scheme defined by \eqref{semi-afc} and \eqref{glow}. The high-order entropy stable IDP scheme of Section \ref{sec:summ} is labeled HO-ES-IDP. The method corresponding to HO-ES-IDP without Step 6 is referred to as HO-ES. We use this version to show the effect of deactivating the IDP limiter. The significance of individual steps of the HO-ES-IDP algorithm is further illustrated by varying the definition of the target flux $f_{ij}^e$ for a two-dimensional test problem with a nonconvex flux function (the so-called KPP problem \cite{kpp}).
In all numerical examples, we discretize in time using the third-order
explicit SSP Runge-Kutta method with three stages \cite{ssprev}.
Unless otherwise stated, we use structured triangular meshes. 
All computations are performed using Proteus
(https://proteustoolkit.org), an open-source Python toolkit for numerical simulations.

\subsection{One-dimensional advection}
The first problem that we consider in this work is the one-dimensional linear advection equation
\beq\label{linad1d}
\pd{u}{t}+a\pd{u}{x}=0\quad\mbox{in}\quad\Omega=(0,1)
\eeq
with the constant velocity $a=1$.
The smooth initial condition is given by 
\beq
u_0(x)=\cos\left(2\pi(x-0.5)\right).
\eeq
We solve \eqref{linad1d} up to the final time $t=1$ and measure
the numerical errors w.r.t. the $L^1$ norm.
In Table \ref{table:advection}, we show the results of a grid convergence study.
As expected, LO-ES-IDP exhibits first-order convergence behavior, while
second-order convergence is achieved with HO-ES and HO-ES-IDP. 

\begin{table}[!h]\scriptsize
  \begin{center}
    \begin{tabular}{|c||c|c||c|c||c|c||} \cline{1-7}
      \multicolumn{1}{|c||}{} &
      \multicolumn{2}{|c||}{LO-ES-IDP} &
      \multicolumn{2}{|c||}{HO-ES} &
      \multicolumn{2}{|c||}{HO-ES-IDP} \\ \hline
      $N_h$ &
      $\|u_h-u_\text{exact}\|_{L^1}$ & EOC &
      $\|u_h-u_\text{exact}\|_{L^1}$ & EOC  &
      $\|u_h-u_\text{exact}\|_{L^1}$ & EOC \\ \hline
      11  & 3.70E-2 &  --  & 1.21E-2 &  --  & 1.69E-2 & --   \\ \hline
      16  & 2.92E-2 & 0.58 & 6.01E-3 & 1.73 & 8.32E-3 & 1.74 \\ \hline
      21  & 2.40E-2 & 0.68 & 3.73E-3 & 1.65 & 5.15E-3 & 1.66 \\ \hline
      31  & 1.77E-2 & 0.75 & 1.82E-3 & 1.76 & 2.51E-3 & 1.77 \\ \hline
      41  & 1.40E-2 & 0.81 & 1.08E-3 & 1.81 & 1.47E-3 & 1.84 \\ \hline
      61  & 9.84E-3 & 0.86 & 5.06E-4 & 1.86 & 6.93E-4 & 1.86 \\ \hline
      81  & 7.59E-3 & 0.90 & 2.95E-4 & 1.88 & 4.01E-4 & 1.90 \\ \hline
      121 & 5.21E-3 & 0.92 & 1.36E-4 & 1.89 & 1.87E-4 & 1.88 \\ \hline
      161 & 3.97E-3 & 0.94 & 7.86E-5 & 1.91 & 1.08E-4 & 1.90 \\ \hline
      241 & 2.68E-3 & 0.96 & 3.65E-5 & 1.89 & 4.94E-5 & 1.92 \\ \hline
      321 & 2.03E-3 & 0.97 & 2.11E-5 & 1.90 & 2.82E-5 & 1.94 \\ \hline
      481 & 1.36E-3 & 0.98 & 9.69E-6 & 1.91 & 1.28E-5 & 1.95 \\ \hline
    \end{tabular}
    \caption{One-dimensional advection. Grid convergence history for three entropy-stable AFC schemes. 
      \label{table:advection}}
  \end{center}
\end{table}

\subsection{Two-dimensional advection}
\begin{figure}[!h]
  \centering
  \subfloat[zoom of the grid\label{fig:xos_mesh}]
           {
             \begin{tabular}{c}
               \\
               \includegraphics[scale=0.12]{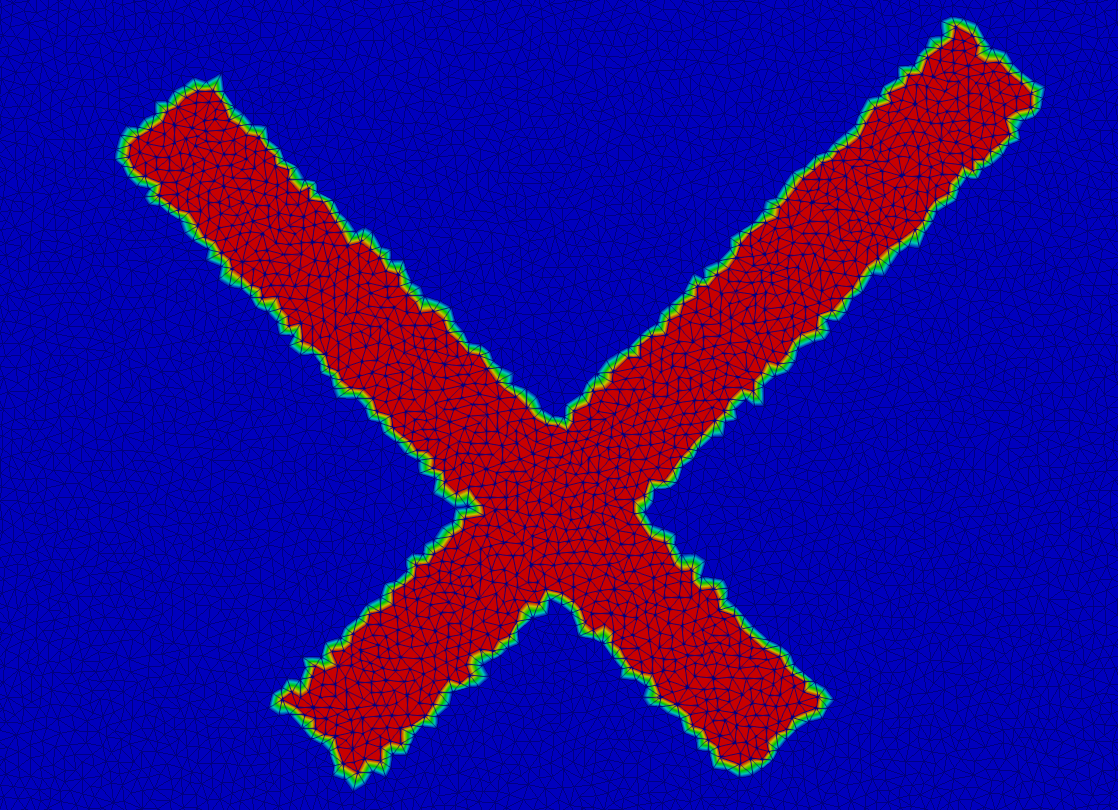}
             \end{tabular}
           }
  \subfloat[initial condition\label{fig:xos_init_cond}]
           {
             \begin{tabular}{c}
               \\
               \includegraphics[scale=0.17]{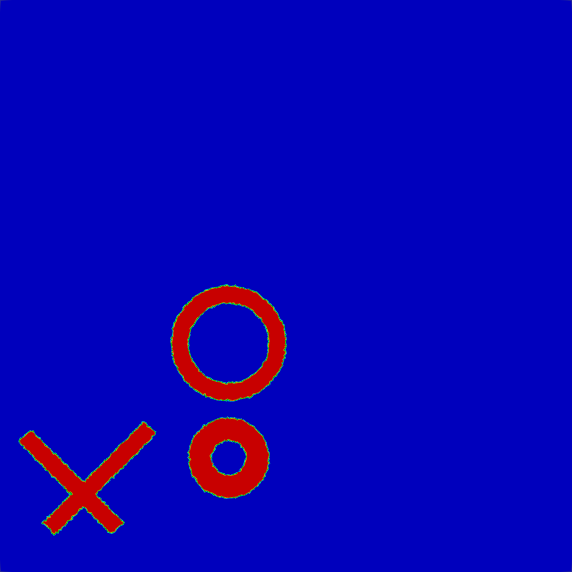} 
             \end{tabular}
           }
    \subfloat[LO-ES-IDP solution\label{fig:low_order}]
             {
               \begin{tabular}{c}
                 {\scriptsize $u^{\max}=0.5034$} \\
                 \includegraphics[scale=0.17]{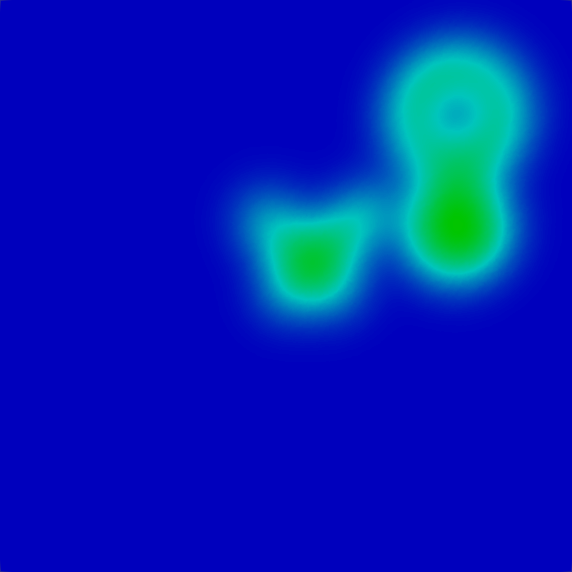}
               \end{tabular}
             }
  \subfloat[HO-ES solution \label{fig:high_order_no_IDP}]
           {
             \begin{tabular}{c}
               {\scriptsize $u^{\max}=0.9862$} \\
               \includegraphics[scale=0.17]{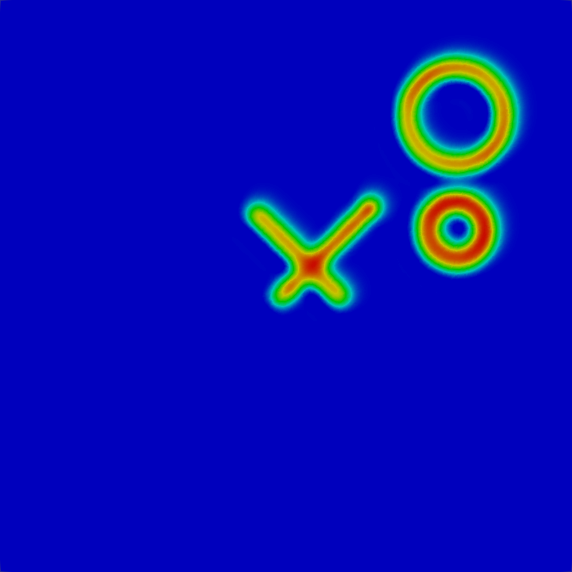}
             \end{tabular}
           }
           \caption{Two-dimensional linear advection problem. Zoom of the unstructured grid, initial
             data, and numerical solutions at $t=4$ obtained with $N_h=99,412$ DoFs.}
\end{figure}

The next example was used in \cite{DG-BFCT} to study the numerical behavior of 
flux-corrected transport algorithms for high-order DG discretizations of the
two-dimensional linear advection problem 
\beq
\pd{u}{t}+{\bf v}\cdot\nabla u=0\quad\mbox{in}\quad\Omega=(0,100)^2
\eeq
with constant velocity ${\bf v}=(10,10)$.
The initial condition, which is shown in Figure \ref{fig:xos_init_cond},
is composed of two rings and a cross. 
The upper ring is centered at $(x,y)=(40,40)$. The radii of its inner and outer
circles are $7$ and $10$, respectively. The center of the lower ring is located
at the point $(x,y)=(40,20)$. The radii of the inner and outer circles are $3$ and $7$,
respectively. The cross occupies the region $r_1\cup r_2\subset\Omega$,
where $r_1=\{x,y\in\Omega ~|~ x\in[7,32], y\in[10,13]\}$
and $r_2=\{x,y\in\Omega ~|~ x\in[14,17], y\in[3,26]\}$,
rotated by $-45^\circ$ around the point $(x,y)=(15.5,11.5)$. 

\begin{figure}[!h]
  \centering
  \begin{tabular}{ccc}
    ${\scriptsize u^{\max}=0.9382}$ &
    ${\scriptsize u^{\max}=0.9947}$ &
    ${\scriptsize u^{\max}=0.9999}$ \\
  \includegraphics[scale=0.2]{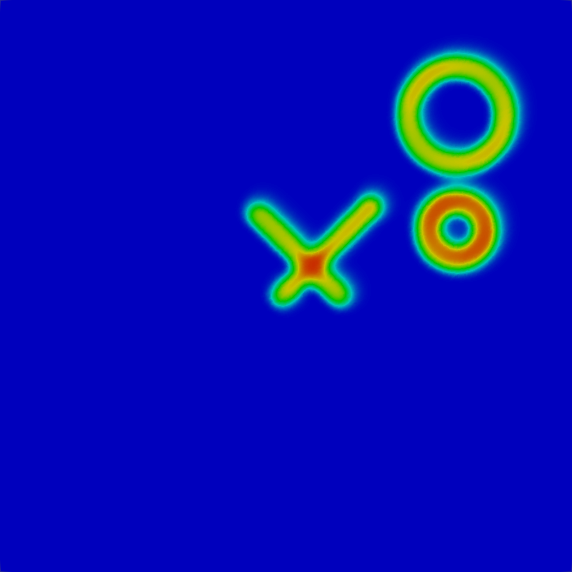} &
  \includegraphics[scale=0.2]{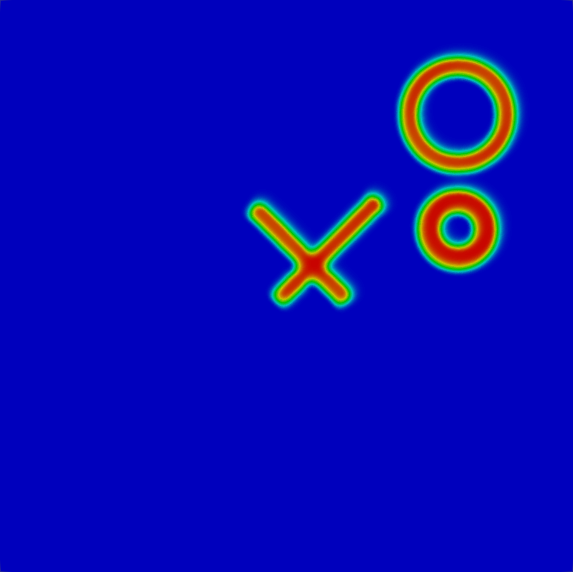} &
  \includegraphics[scale=0.2]{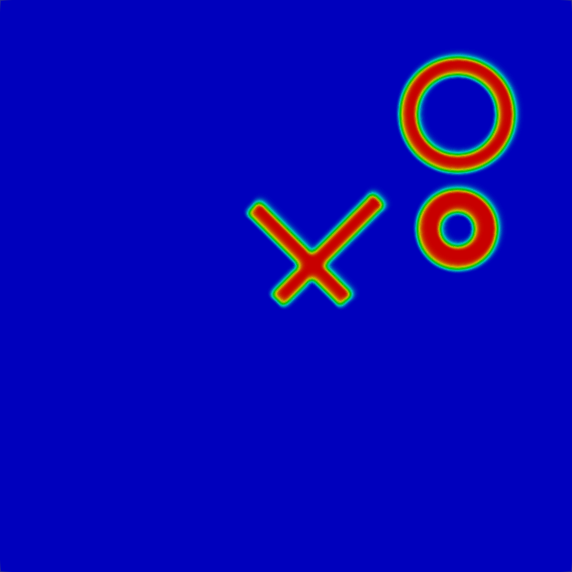} \\
  (a) $N_h=99,412$ & (b) $N_h=395,745$ & (c) $N_h=1,580,651$
  \end{tabular}
  \caption{Two-dimensional linear advection problem. Numerical solutions
    at $t=4$ obtained with HO-ES-IDP on three successively refined
    unstructured grids. 
    \label{fig:xos_mult_refinements}}
\end{figure}

For this problem, we use unstructured grids. In Figure \ref{fig:xos_mesh}, we show a zoom
of one of these grids. Computations are terminated at the final time $t=4$. In Figures
\ref{fig:low_order} and \ref{fig:high_order_no_IDP}, we present the LO-ES-IDP and HO-ES
solutions calculated using $N_h=99,412$ degrees of freedom (DoFs). The higher accuracy
of the latter approximation illustrates the need for antidiffusive corrections of the
LLF flux. No significant undershoots or overshoots are generated by HO-ES in this example.
The results obtained with the HO-ES-IDP scheme on three successively refined meshes are
shown in Figure~\ref{fig:xos_mult_refinements}. The advected discontinuities are resolved
in a crisp and nonoscillatory manner, especially on the finest mesh.

\subsection{KPP problem}
The KPP problem \cite{Guermond2016,Guermond2017,kpp} is a
challenging nonlinear test for verification of entropy stability
properties. We use this problem to 
test different components of the method that we propose. 
In this series of 2D experiments,
we solve equation \eqref{ibvp-pde} with the nonlinear and
nonconvex flux function
\beq
\mathbf{f}(u)=(\sin(u),\cos(u))
\eeq
 in the computational domain
$\Omega=(-2,2)\times(-2.5,1.5)$ using the initial
condition
\beq
u_0(x,y)=\begin{cases}
\frac{14\pi}{4} & \mbox{if}\quad \sqrt{x^2+y^2}\le 1,\\
\frac{\pi}{4} & \mbox{otherwise}.
\end{cases}
\eeq
A simple (but rather pessimistic) upper bound for the guaranteed
maximum speed (GMS) is $\lambda=1$. More accurate GMS estimates
can be found in \cite{Guermond2017}. The exact solution exhibits
a two-dimensional rotating wave structure, which is difficult to
capture in numerical simulations using high-order
methods. The main challenge of this test is to
prevent possible convergence to wrong weak solutions.

\begin{figure}[!h]
  \centering
  \subfloat[\label{fig:kpp_low_order}]{\includegraphics[scale=0.17]{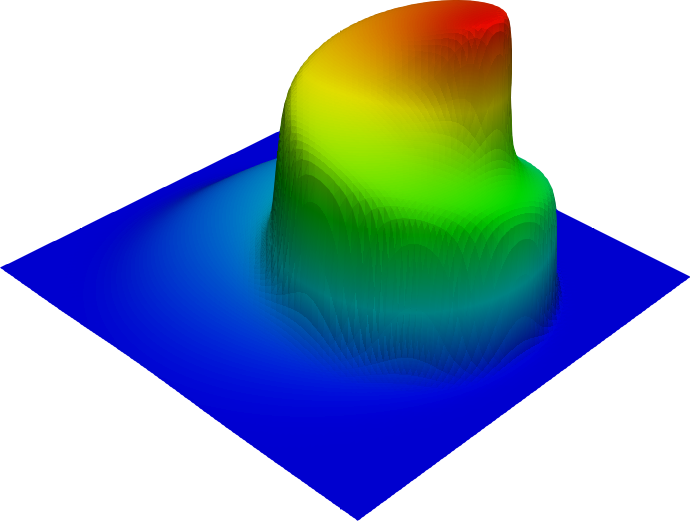}}
  \subfloat[\label{fig:kpp_galerkin}]{\includegraphics[scale=0.17]{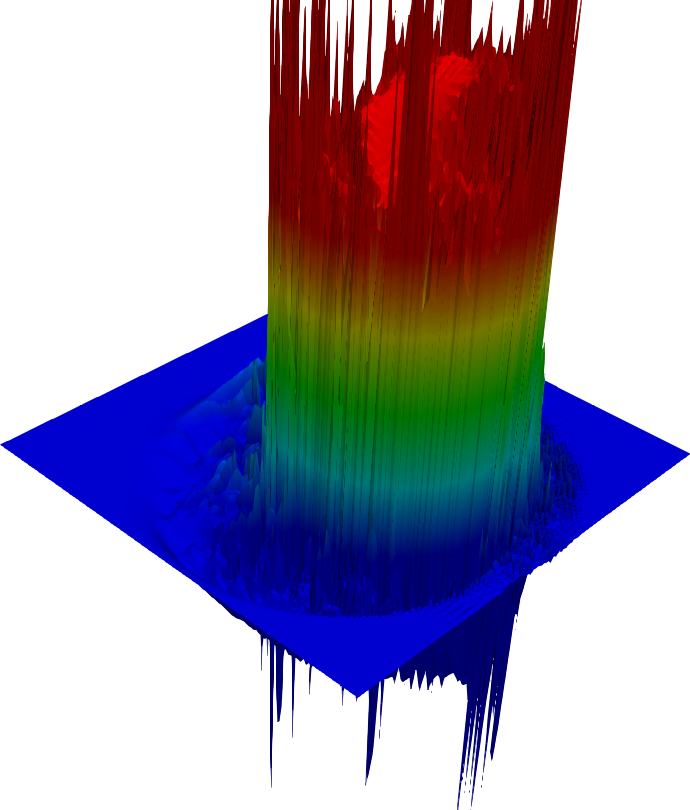}}
  \subfloat[\label{fig:kpp_no_ev_with_IDP}]{\includegraphics[scale=0.17]{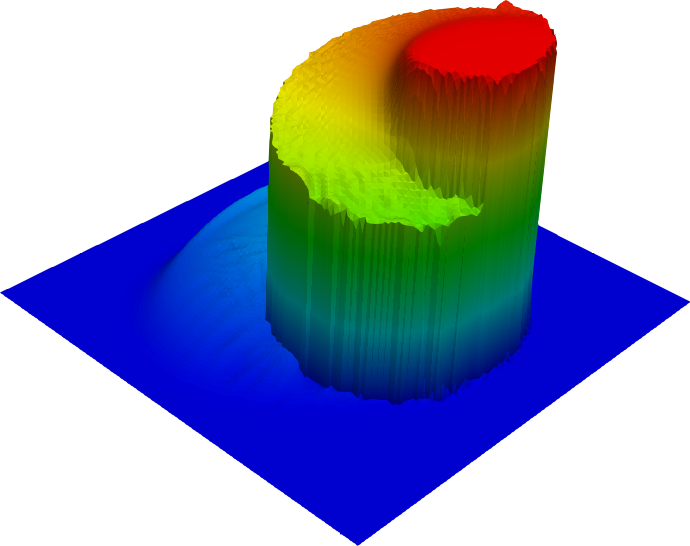}}
  \subfloat[\label{fig:kpp_no_ev_with_MC_IDP}]{\includegraphics[scale=0.17]{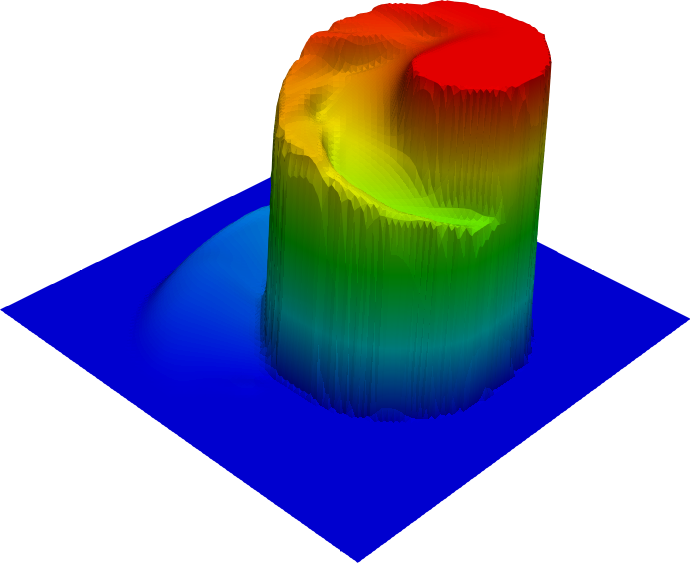}}
           \caption{Numerical solutions of the KPP problem at $t=1$ obtained using $N_h=129^2$ DoFs
             and (a) LO-ES-IDP, (b) unconstrained lumped-mass Galerkin method corresponding to
             $f_{ij}^{e,**}=d_{ij}^{e,\max}(u_i-u_j)=f_{ij}^e$, (c) IDP limiter 
             \eqref{fij_lim} for $f_{ij}^e=d_{ij}^{e,\max}(u_i-u_j)$, and (d) IDP limiter
             \eqref{fij_lim} for $f_{ij}^e=m_{ij}^e(\dot u_i-\dot u_j)+d_{ij}^{e,\max}(u_i-u_j)$.}
\end{figure}

Numerical solutions are evolved up to the final time $t=1$.
To test individual components of the AFC scheme summarized in
Section \ref{sec:summ}, we vary the definition of the
target flux $f_{ij}^e$ and/or the way in which it is limited to produce a
constrained flux
$f_{ij}^{e,**}$. In Figure \ref{fig:kpp_low_order}, we present the LO-ES-IDP solution
($f_{ij}^{e,**}=0=f_{ij}^e$). It is highly dissipative but provides a correct
qualitative description of the rotating wave structure. The solution
displayed in Figure \ref{fig:kpp_galerkin} was calculated
using $f_{ij}^{e,**}=d_{ij}^e(u_i-u_j)=f_{ij}^e$.
This lumped-mass Galerkin approximation
is highly oscillatory and exhibits an entropy-violating
merger of two shocks. The
solutions shown in Figures~\ref{fig:kpp_no_ev_with_IDP} and
\ref{fig:kpp_no_ev_with_MC_IDP}
were obtained using the IDP-limited counterparts $f_{ij}^{e,**}=f_{ij}^{e,*}$
of $f_{ij}^{e}=d_{ij}^e(u_i-u_j)$ and $f_{ij}^{e}=m_{ij}^e(\dot u_i-\dot u_j)
+d_{ij}^e(u_i-u_j)$, respectively. As reported in \cite{convex,CL-diss},
the latter definition of the target flux introduces high-order
background dissipation. However, neither the activation of the
IDP limiter nor the inclusion of $m_{ij}^e(\dot u_i-\dot u_j)$
provides enough entropy dissipation
to capture the twisted shocks correctly in the KPP test. This
unsatisfactory state of affairs
illustrates the need for entropy stabilization and confirms 
the findings of Guermond et al. \cite{Guermond2016} who noticed that
IDP limiting alone does not guarantee convergence to entropy solutions.

\begin{figure}[!h]
  \centering
  \begin{tabular}{ccc}
    \includegraphics[scale=0.22]{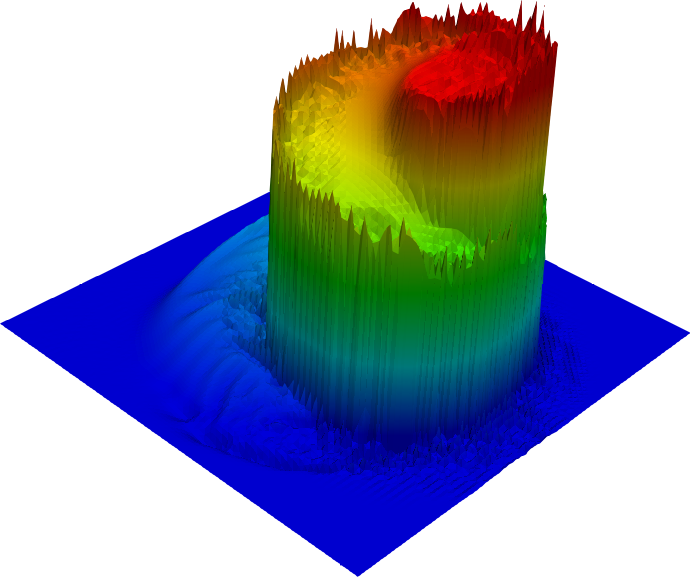} &
    \includegraphics[scale=0.22]{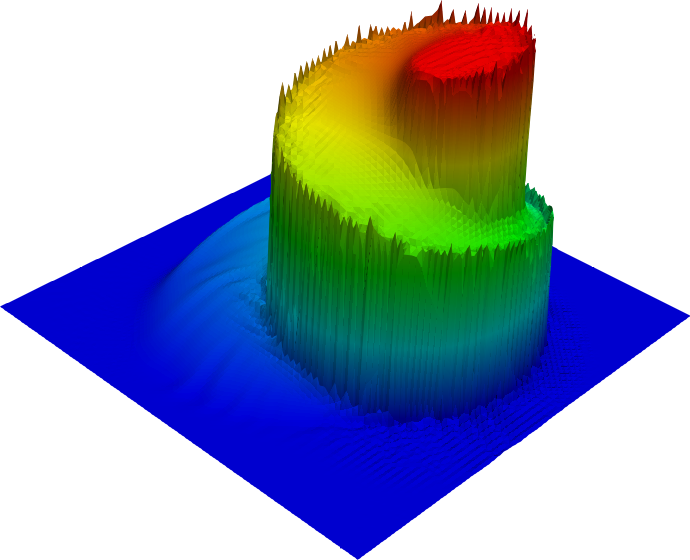} &
    \includegraphics[scale=0.22]{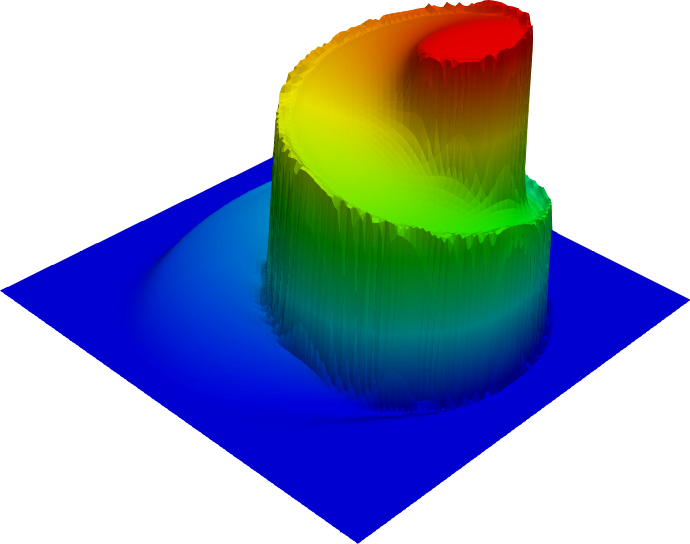}
    \\ & & \\
    \includegraphics[scale=0.22]{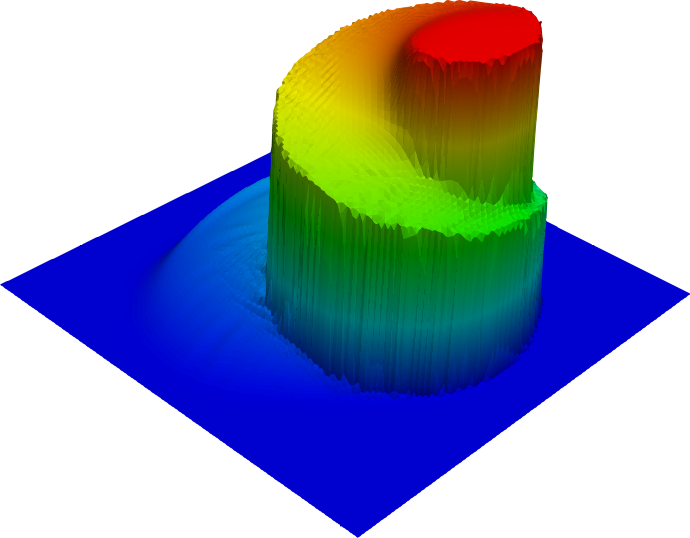} &
    \includegraphics[scale=0.22]{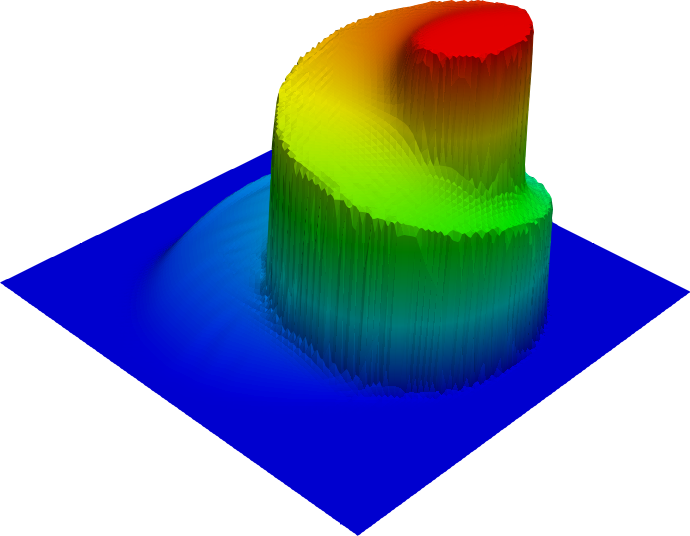} &
    \includegraphics[scale=0.22]{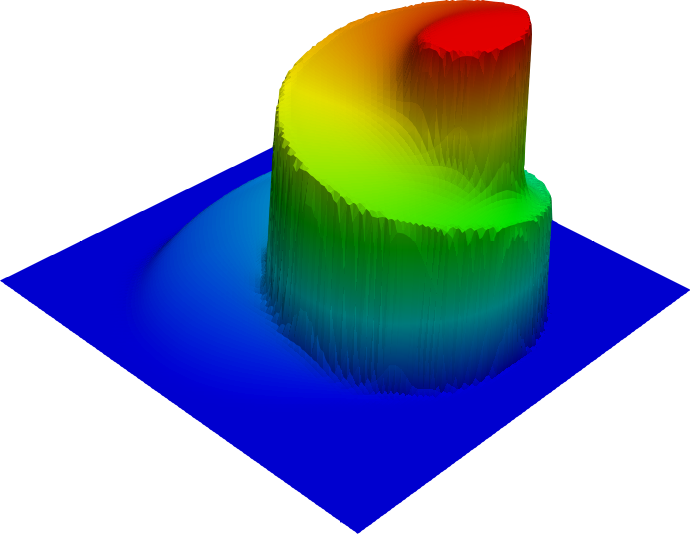}
    \\
    (a) &
    (b) &
    (c)
  \end{tabular}
  \caption{Numerical solutions of the KPP problem  obtained at $t=1$
    using $N_h=129^2$ DoFs. The diagrams of the first and second rows show the
    results produced by three entropy stable high-order schemes without and with
    activation of the IDP flux limiter, respectively. The target fluxes are defined
    by (a) $f_{ij}^e=(d_{ij}^{e,\min}-d_{ij}^{e,\max})(u_j^n-u_i^n)$ for the diagrams of
    the first column, (b) $f_{ij}^e=(d_{ij}^{e,\min}-d_{ij}^{e,\max})(u_j^n-u_i^n)+\nu_{ij}^e(v_j^n-v_i^n)$
    for the diagrams of the second column, and (c)
    $f_{ij}^e=m_{ij}^e(\dot u_i-\dot u_j)+(d_{ij}^{e,\min}+d_{ij}^{e,\max})(u_j-u_i)+\nu_{ij}^e(v_j-v_i)$
    for the diagrams of the third column.
  }\label{fig:kpp_entropy_stable}
\end{figure}

\begin{figure}[!h]
  \centering
  \subfloat[$N_h=257^2$]{\includegraphics[scale=0.22]{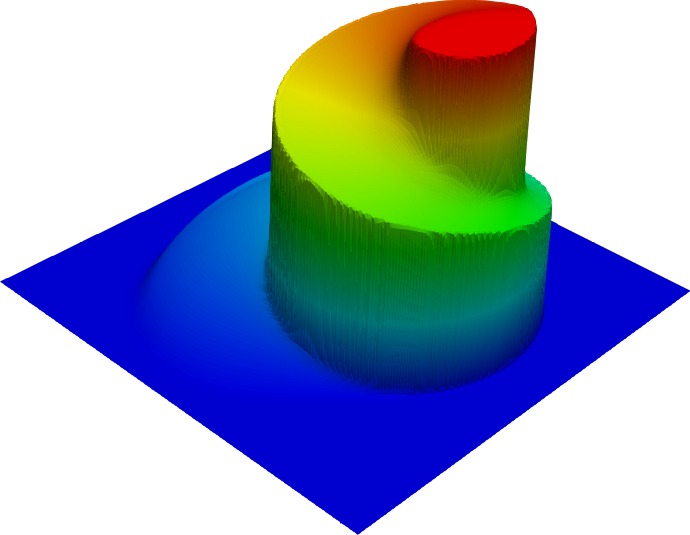}}
  \subfloat[$N_h=513^2$]{\includegraphics[scale=0.22]{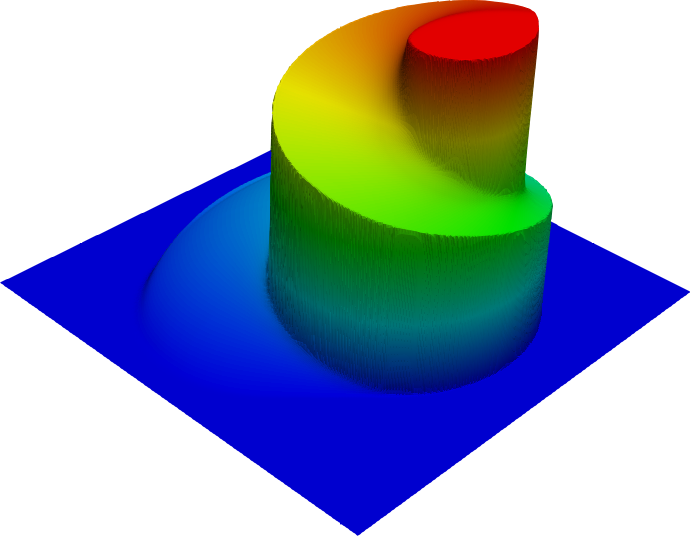}}
  \subfloat[$N_h=1025^2$]{\includegraphics[scale=0.22]{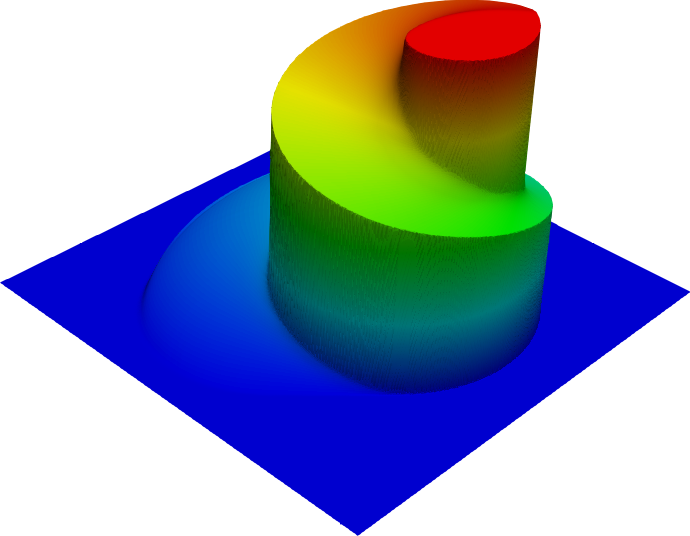}}
  \caption{Numerical solutions of the KPP problem at $t=1$ obtained using HO-ES-IDP on three meshes.}\label{fig:kpp_mult_refinements}
\end{figure}

In Figure \ref{fig:kpp_entropy_stable}, we present numerical solutions produced
by three entropy stable high-order methods using $N_h=129^2$. Each column corresponds to a different
definition of the target flux $f_{ij}^e$. The diagrams of the first row were
calculated without invoking the IDP flux limiter \eqref{fij_lim}. Therefore, the
results exhibit undershoots/overshoots. The target fluxes of the three schemes
under investigation are listed in the caption. The use of
$f_{ij}^e=(d_{ij}^{e,\min}-d_{ij}^{e,\max})(u_j-u_i)$ produces a barely
entropy stable approximation which requires additional entropy fixes
to  keep the shocks separated if no IDP constraints are imposed. The inclusion
of $\nu_{ij}^e(v_j^n-v_i^n)$ and $m_{ij}^e(\dot u_i-\dot u_j)$ leads to AFC schemes
that reproduce the rotating wave structure correctly even without IDP limiting.
The method employed in the last diagram of the second row is HO-ES-IDP. In Figure
\ref{fig:kpp_mult_refinements}, we show the HO-ES-IDP results for finer
meshes. As the number of degrees of freedom is increased, the method
converges to a bound-preserving entropy solution.

\subsection{Inviscid Burgers equation}

Let us now consider the two-dimensional
inviscid Burgers equation \cite{GuermondNazarov2014,convex}
\beq
\pd{u}{t}+\nabla\cdot\left(\mathbf{v}\frac{u^2}{2}\right)=0\qquad
\mbox{in}\ \Omega=(0,1)^2, 
\eeq
where $\mathbf{v}=(1,1)$ is a constant vector. The
piecewise-constant initial data is given by
\beq
u_0(x,y)=\begin{cases}
-0.2 & \mbox{if}\quad x < 0.5\ \land y >0.5,\\
-1.0 & \mbox{if}\quad x > 0.5\ \land y >0.5,\\
\phantom{-}0.5 & \mbox{if}\quad x < 0.5\ \land y <0.5,\\
\phantom{-}0.8 & \mbox{if}\quad x > 0.5\ \land y <0.5.
\end{cases}
\eeq
The inflow boundary conditions are defined using the exact solution of the
pure initial value problem in~$\R^2$.  This solution can be found in
\cite{GuermondNazarov2014} and stays in the invariant set
$\mathcal G=[-1.0,0.8]$.

\begin{figure}[!h]
  \centering
 \subfloat[$N_h=33^2$]{\includegraphics[scale=0.275]{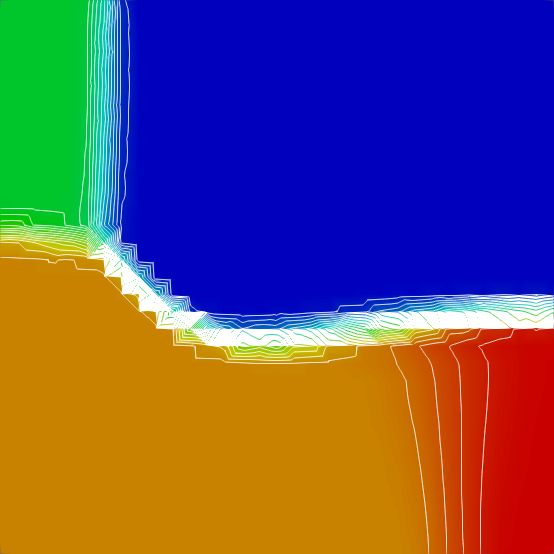}}
 \subfloat[$N_h=129^2$]{ \includegraphics[scale=0.275]{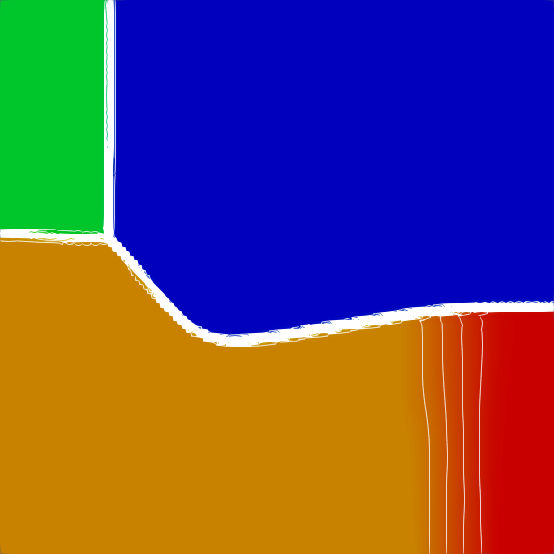}}
 \subfloat[$N_h=513^2$]{\includegraphics[scale=0.275]{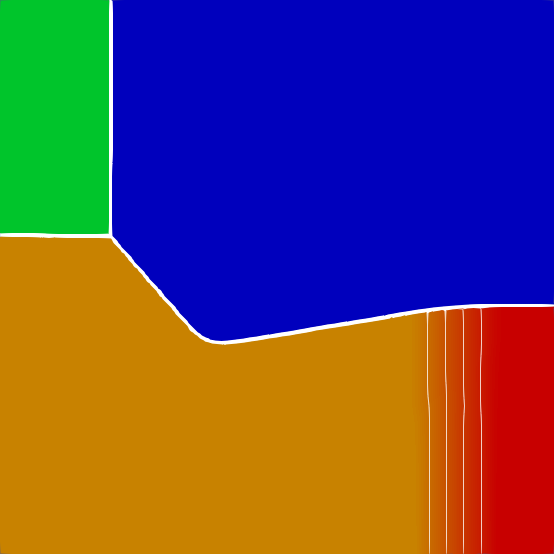}} 
 \caption{Inviscid Burgers equation in two dimensions. Numerical solutions
    at $t=0.5$ obtained with HO-ES-IDP on three meshes. In each diagram, we
   plot 30 contour lines corresponding to a uniform subdivision of
   $\mathcal G=[-1.0,0.8]$.
  }\label{fig:burgers}
\end{figure}

The final time for computation of numerical solutions is $t=0.5$. In
Table \ref{table:burgers}, we show the results of a grid convergence
study for LO-ES-IDP, HO-ES, and HO-ES-IDP. In this example, the low-order
LLF scheme performs remarkably well at self-steepening shocks but the
resolution of rarefactions is not as accurate as in the case of the
high-order entropy stable method without and with IDP limiting. The
HO-ES-IDP results calculated on three successively refined meshes are
displayed in Figure \ref{fig:burgers}. 

\begin{table}[!h]\scriptsize
  \begin{center}
    \begin{tabular}{|c||c|c||c|c||c|c||} \cline{1-7}
      \multicolumn{1}{|c||}{} &
      \multicolumn{2}{|c||}{LO-ES-IDP} &
      \multicolumn{2}{|c||}{HO-ES} &
      \multicolumn{2}{|c||}{HO-ES-IDP} \\ \hline
      $N_h$ &
      $\|u_h-u_\text{exact}\|_{L^1}$ & EOC &
      $\|u_h-u_\text{exact}\|_{L^1}$ & EOC  &
      $\|u_h-u_\text{exact}\|_{L^1}$ & EOC \\ \hline
      $33^2$   & 7.63E-2 &  --  & 4.02E-2 &--& 3.93E-2 & --   \\ \hline
      $65^2$   & 4.49E-2 & 0.76 & 2.12E-2 & 0.92 & 2.09E-2 & 0.91  \\ \hline
      $129^2$  & 2.51E-2 & 0.85 & 1.11E-2 & 0.95 & 1.10E-2 & 0.95 \\ \hline
      $257^2$  & 1.37E-2 & 0.85 & 5.57E-3 & 0.97 & 5.62E-3 & 0.94 \\ \hline
      $513^2$  & 7.31E-3 & 0.90 & 2.80E-3 & 0.99 & 2.83E-3 & 0.98 \\ \hline
    \end{tabular}
    \caption{Inviscid Burgers equation in two dimensions. 
      Grid convergence history for three entropy stable AFC schemes. 
      \label{table:burgers}}
  \end{center}
\end{table}

\subsection{Buckley-Leverett equation}
\begin{figure}[!h]
  \centering
  \begin{tabular}{ccc}
    ${\scriptsize u^{\max}=0.9469}$ &
    ${\scriptsize u^{\max}=1.0002}$ &
    ${\scriptsize u^{\max}=0.9923}$ \\
    \includegraphics[scale=0.25]{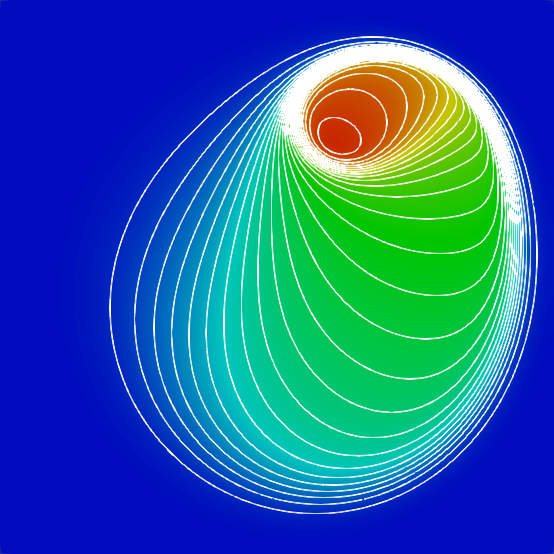} &
    \includegraphics[scale=0.25]{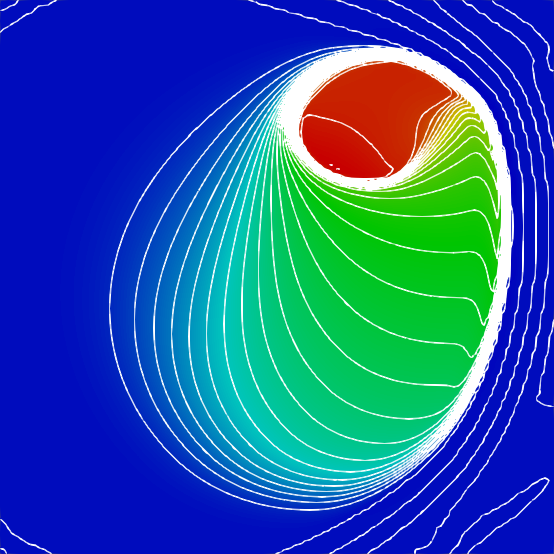} &
    \includegraphics[scale=0.25]{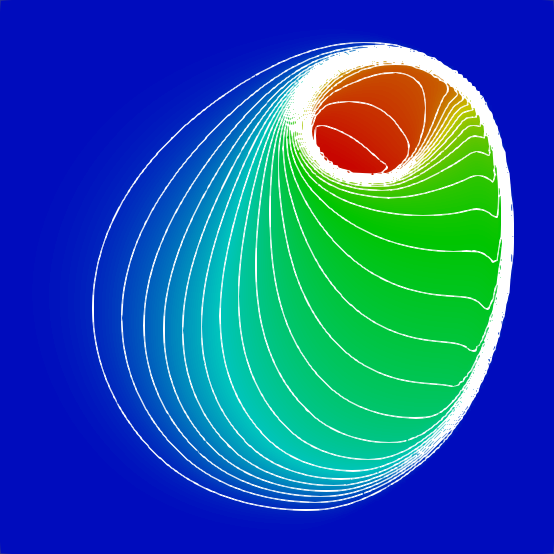} \\
    (a) LO-ES-IDP & (b) HO-ES & (c) HO-ES-IDP
  \end{tabular}
  \caption{Two-dimensional Buckley-Leverett problem. Numerical solutions at $t=0.5$
    obtained with three entropy-stable AFC schemes using $\mathbb{Q}_1$ finite elements
    and $N_h=129^2$ DoFs. In each diagram, we plot 30 contour lines corresponding
    to a uniform subdivision of $\mathcal G=[0,1]$.
    \label{fig:bl_methods}}
\end{figure}

\begin{figure}[!h]
  \centering
  \begin{tabular}{ccc}
    ${\scriptsize u^{\max}=1}$ &
    ${\scriptsize u^{\max}=1}$ &
    ${\scriptsize u^{\max}=1}$ \\
    \includegraphics[scale=0.25]{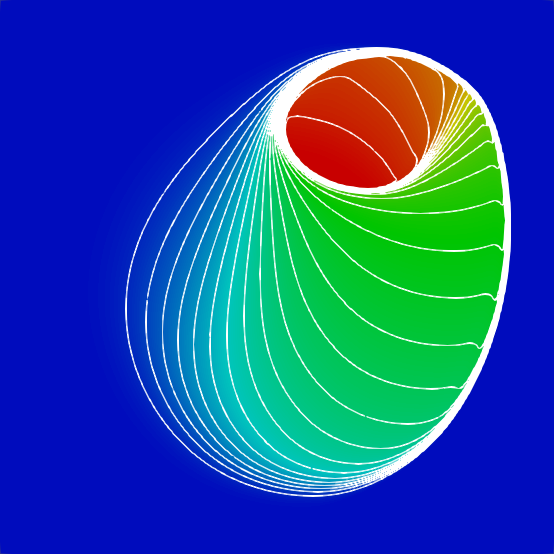} &
    \includegraphics[scale=0.25]{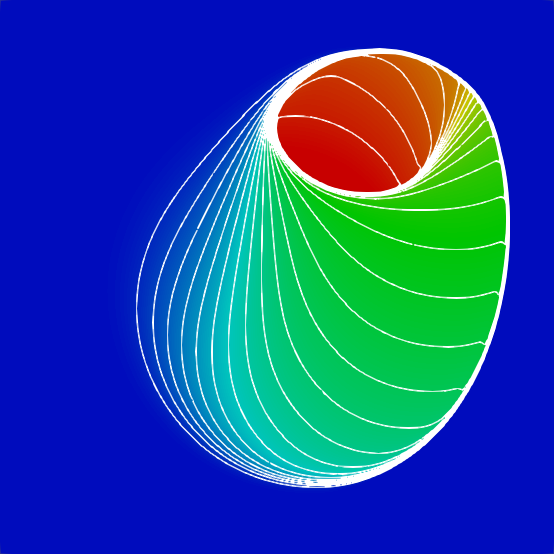} &
    \includegraphics[scale=0.25]{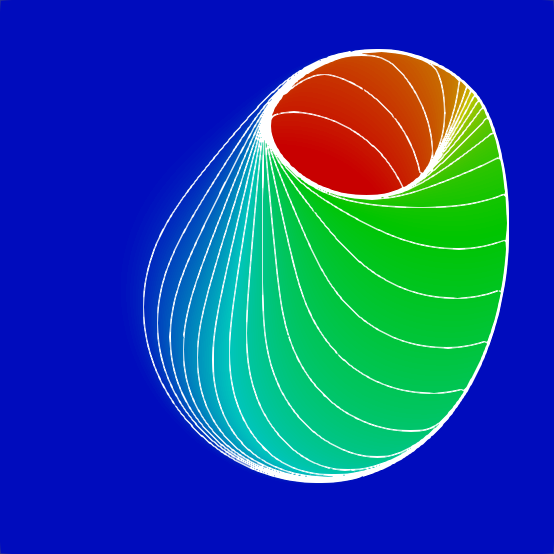} 
  \end{tabular}
  \caption{
    Two-dimensional Buckley-Leverett problem. Numerical solutions at $t=0.5$ obtained
    with HO-ES-IDP on three meshes. In each diagram, we plot 30 contour lines corresponding
    to a uniform subdivision of $\mathcal G=[0,1]$.
    \label{fig:bl_refinements}}
\end{figure}

In the last numerical experiment, we consider the two-dimensional
Buckley-Leverett equation. The nonconvex
flux function of the nonlinear conservation law to be solved
is  \cite{christov2008new}
\beq
\mathbf{f}(u)=
\frac{u^2}{u^2+(1-u)^2}
\begin{bmatrix}
1 \\ 1-5(1-u)^2
\end{bmatrix}.
\eeq
The computational domain is $\Omega=(-1.5,1.5)^2$. The initial condition is given by
\beq
u_0(x,y)=\begin{cases}
1 & \mbox{ if } x^2+y^2<0.5, \\
0 & \mbox{ otherwise}.
\end{cases}
\eeq
An upper bound for the fastest wave speed can be found in \cite{christov2008new}. 
Similarly to the KPP problem, the solution exhibits a rotating wave structure.
In Figure \ref{fig:bl_methods}, we show the entropy stable AFC
approximations at the final time
$t=0.5$. Note that small oscillations are present if the IDP
flux limiter is not applied. The effect of mesh refinement on the
accuracy of HO-ES-IDP is illustrated by the snapshots presented in Figure
\ref{fig:bl_refinements}. In all experiments for this test problem
we use bilinear finite elements. 

\section{Conclusions}

We have shown that algebraic flux correction schemes can be configured
to satisfy discrete entropy inequalities in addition to 
discrete maximum
principles. The new inequality-constrained
stabilization techniques modify the residual of the
semi-discrete Galerkin
scheme in a way which ensures entropy stability while preserving
all other important properties (conservation, preservation of
local bounds, low levels of numerical diffusion). The proposed
methodology was presented in the context of continuous
Galerkin methods. Further developments will focus on the DG version
\cite{abgrall,ranocha,chen,pazner}, extensions to high-order finite
elements \cite{DG-BFCT,convex2}, and design of entropy stability
preserving time integrators \cite{relax}.

\medskip
\paragraph{\bf Acknowledgments}

The work of Dmitri Kuzmin was supported by the German Research Association (DFG) under grant KU 1530/23-1. The authors would like to thank Hennes Hajduk (TU Dortmund University) for careful proofreading of the manuscript and helpful feedback.

%\section*{References}

\end{document}